\documentclass[a4paper,10pt]{amsart}
\usepackage{amsmath}
\usepackage{amssymb}
\usepackage{color}
\usepackage{a4wide}

\theoremstyle{plain}
\newtheorem{theorem}{Theorem}
\newtheorem{corollary}{Corollary}
\newtheorem{lemma}{Lemma}

\newtheorem{remark}{Remark}

\begin{document}

\title[Functional a posteriori error estimates for time-periodic parabolic problems]{
Functional a posteriori error estimates for parabolic time-periodic boundary value problems
}

\author{Ulrich Langer}
\author{Sergey Repin}
\author{Monika Wolfmayr}

\address[U. Langer]{Institute of Computational Mathematics,
Johannes Kepler University Linz,
Altenbergerstra\ss e 69, 4040 Linz, Austria}
\email{ulanger@numa.uni-linz.ac.at}

\address[S. Repin]{V. A. Steklov Institute of Mathematics in St. Petersburg,
Fontanka 27, 191011, St. Petersburg, Russia,
and University of Jyv\"{a}skyl\"{a}, Finland}
\email{repin@pdmi.ras.ru}

\address[M. Wolfmayr]{Johann Radon Institute for Computational and Applied Mathematics,
Altenbergerstra\ss e 69, 4040 Linz, Austria}
\email{monika.wolfmayr@ricam.oeaw.ac.at}

\begin{abstract}
The paper is concerned with
parabolic time-periodic boundary value problems which
are of theoretical interest and arise in different practical
applications. The multiharmonic finite element method is 
well adapted to this class of parabolic problems.
We study properties of multiharmonic approximations and derive
guaranteed and fully computable bounds of approximation errors.
For this purpose, we use the functional a posteriori error estimation
techniques earlier
introduced by S. Repin.
Numerical tests confirm the efficiency of the 
a posteriori error bounds derived.
\end{abstract}

\maketitle

\section{Introduction}
\label{Sec1:Intro}

Initial-boundary value problems for parabolic equations describe many
quite different physical phenomena such as
heat conduction, diffusion, chemical reactions,
biological processes, and transient electromagnetical fields. 
The numerical simulation of these
phenomena is usually based on time-integration methods together with a suitable space
discretization, see, e.g., the well-known monograph
\cite{LRW:Thomee:2006} and the references therein.
In many practically interesting cases, for instance,
in electromagnetics and chemistry,
the processes are time-periodic, see, e.g.,
\cite{LRW:AbbeloosDiehlHinzeVandewalle:2011}.
In this case, the initial condition must be replaced by the time-periodicity condition.
Standard time-integration methods may be less efficient then methods based on approximations
in terms of Fourier series.
This paper deals with this type of
approximations. 
In fact, it is devoted to the a posteriori error analysis of parabolic
time-periodic boundary value problems in connection with their multiharmonic finite element
discretization. More precisely, all functions are expanded into Fourier series, approximations
are presented by truncated series and the Fourier coefficients
are approximated by the finite element method (FEM). This so-called
multiharmonic FEM (MhFEM) or harmonic-balanced FEM was successfully used for the simulation
of electromagnetic devices described by nonlinear eddy current problems with harmonic
excitations, see, e.g., \cite{LRW:YamadaBessho:1988,
LRW:BachingerLangerSchoeberl:2005,
LRW:BachingerLangerSchoeberl:2006,
LRW:CopelandLanger:2010} and the references therein. Later, this
discretization technique has been applied to linear time-periodic parabolic boundary value and optimal
control problems \cite{LRW:KollmannKolmbauer:2011,
LRW:KollmannKolmbauerLangerWolfmayrZulehner:2013, LRW:KrendlSimonciniZulehner:2012,
LRW:LangerWolfmayr:2013, LRW:Wolfmayr:2014} and to linear time-periodic eddy current problems and
the corresponding optimal control problems
\cite{LRW:Kolmbauer:2012c, LRW:KolmbauerLanger:2012, LRW:KolmbauerLanger:2013}.
In this framework,
we deduce a posteriori error estimates which provide guaranteed and fully computable upper bounds
(majorants) of the respective errors.
To the best of our knowledge these estimates are new.
Our approach is based on the works of Repin, see, e.g., the papers on
parabolic problems \cite{LRW:Repin:2002,
LRW:GaevskayaRepin:2005} as well as on optimal control problems
\cite{LRW:GaevskayaHoppeRepin:2006,
LRW:GaevskayaHoppeRepin:2007},
the books \cite{LRW:Repin:2008,
LRW:MaliNeittaanmaekiRepin:2014},
and the references therein. In particular, our a posteriori error analysis
uses the techniques
close to the one suggested in \cite{LRW:Repin:2002}, but the analysis contains essential
changes. In the MhFEM setting, we are able to establish inf-sup and sup-sup conditions from
which we deduce existence and uniqueness of the solution to the parabolic time-periodic
problems by applying the theorem of Babu\v{s}ka and Aziz. Then, we deduce the a posteriori estimates,
which are very valuable for the evaluation of quality of the multiharmonic finite element
solution because they can judge on the quality
of approximation for any particular harmonic. This is highly important
because for linear time-periodic parabolic problems, the computations of the
Fourier coefficients corresponding to every single mode $k = 0, 1, \dots$ are
decoupled. Hence, we can use different meshes independently generated by adaptive finite element
approximations to the Fourier coefficients for different modes. Then, by prescribing certain
bounds, we can finally filter out the Fourier coefficients, which are important for the
numerical solution of the problem. Altogether, such an adaptive multiharmonic finite
element method (AMhFEM) yields complete adaptivity in space and time.
This work is a starting point for the construction of this
AMhFEM, which utilizes the above principles. However, in this work
we are not focused on mesh adaptation issues.
This will be the subject of a separate paper.
Our goal is to provide a detailed a posteriori error
analysis of a parabolic time-periodic boundary value problem in the context of the MhFEM 
leading to guaranteed, computable upper bounds with efficiency indices close to one.

The paper is organized as follows. In Section~\ref{Sec2:ParabolicTimePerBVP},
we discuss a space-time variational formulation 
for parabolic time-periodic boundary value problems
that forms the basis of the MhFEM considered in Section~\ref{Sec3:MhFEApprox}.
Section~\ref{Sec4:FunctionalAPostErrorEstimates}
is devoted to the derivation of functional type a posteriori error
estimates adapted to problems in question.
Finally, in Section~\ref{Sec5:NumericalResults}, we discuss some implementation issues and present first
numerical results.

\section{A parabolic time-periodic boundary value problem}
\label{Sec2:ParabolicTimePerBVP}

Let $Q_T := \Omega \times (0,T)$ denote the space-time cylinder and $\Sigma_T := \Gamma \times (0,T)$
its mantle boundary, where $\Omega \subset \mathbb{R}^d$, $d \in \{1,2,3\}$, is a bounded Lipschitz domain
with the boundary $\Gamma$,
and $(0,T)$ is a given time interval.
The following parabolic time-periodic boundary value problem is considered:
Find $u$ such that
\begin{align}
 \label{problem:parPeriodicBVP:PDE}
 \sigma(\boldsymbol{x}) \, \partial_t u(\boldsymbol{x},t)
	- \text{div} \, (\nu(\boldsymbol{x}) \, \nabla u(\boldsymbol{x},t))
	&= f(\boldsymbol{x},t) \hspace{1cm} &(\boldsymbol{x},t) \in Q_T, \\
 \label{problem:parPeriodicBVP:BC}
 u(\boldsymbol{x},t) &= 0 \hspace{1cm} &(\boldsymbol{x},t) \in \Sigma_T, \\
 \label{problem:parPeriodicBVP:PerCond}
 u(\boldsymbol{x},0) &= u(\boldsymbol{x},T) \hspace{1cm} &\boldsymbol{x} \in \overline{\Omega},
\end{align}
where $f(\boldsymbol{x},t)$ is a given function in $L^2(Q_T)$, and
$\sigma(\cdot)$ and $\nu(\cdot)$ satisfy the assumptions 
\begin{align}
 \label{assumptions:sigmaNu:sigmaStrictlyPositive}
 0 < \underline{\sigma} \leq \sigma(\boldsymbol{x}) \leq \overline{\sigma}, \qquad
 0 < \underline{\nu} \leq \nu(\boldsymbol{x}) \leq \overline{\nu}, \qquad
 \boldsymbol{x} \in \Omega.
\end{align}
In order to study the parabolic time-periodic boundary value problem
(\ref{problem:parPeriodicBVP:PDE})-(\ref{problem:parPeriodicBVP:PerCond}),
we will derive space-time variational formulations in 
Sobolev spaces of functions in the space-time cylinder $Q_T$
using the approach similar to that used by Ladyzhenskaya et al.,
see \cite{LRW:Ladyzhenskaya:1973, LRW:LadyzhenskayaSolonnikovUralceva:1968}.
Let the Sobolev spaces
$H^{1,0}(Q_T) = \{u \in L^2(Q_T) : \nabla u \in [L^2(Q_T)]^d \}$
and $H^{1,1}(Q_T) = \{u \in L^2(Q_T) : \nabla u \in [L^2(Q_T)]^d, \partial_t u \in L^2(Q_T) \}$
be equipped with the norms
\begin{align*}
 \|u\|_{H^{1,0}(Q_T)} &:= \left(\int_{Q_T}
			\left(u(\boldsymbol{x},t)^2 + |\nabla u(\boldsymbol{x},t)|^2 \right)
			\, d\boldsymbol{x} \, dt \right)^{1/2}\; \mbox{and}\\
 \|u\|_{H^{1,1}(Q_T)} &:= \left(\int_{Q_T} \left(u(\boldsymbol{x},t)^2 + |\nabla u(\boldsymbol{x},t)|^2
			+ |\partial_t u(\boldsymbol{x},t)|^2 \right) \, d\boldsymbol{x} \, dt \right)^{1/2},
\end{align*}
respectively, where 
$\nabla = \nabla_{\boldsymbol{x}}$ and $\partial_t$ denote the generalized
derivatives with respect to $x$ and $t$.
The Sobolev space
$H^{0,1}(Q_T) = \{u \in L^2(Q_T) : \partial_t u \in L^2(Q_T) \}$ is defined
analogously. 
Furthermore, the boundary and time-periodicity conditions
are included by defining the Sobolev spaces
\begin{align*}
 H^{1,0}_0(Q_T) &= \{u \in H^{1,0}(Q_T) : u = 0 \text{ on } \Sigma_T \}, \\
 H^{1,1}_0(Q_T) &= \{u \in H^{1,1}(Q_T) : u = 0 \text{ on } \Sigma_T \}, \\
 H^{0,1}_{per}(Q_T) &= \{u \in H^{0,1}(Q_T) : u(\boldsymbol{x},0) = u(\boldsymbol{x},T) \text{ for almost all }
 \boldsymbol{x} \in \Omega \}, \\
 H^{1,1}_{per}(Q_T) &= \{u \in H^{1,1}(Q_T) : u(\boldsymbol{x},0) = u(\boldsymbol{x},T) \text{ for almost all }
 \boldsymbol{x} \in \Omega \}, \\
 H^{1,1}_{0,per}(Q_T) &= \{u \in H^{1,1}_0(Q_T) : u(\boldsymbol{x},0) = u(\boldsymbol{x},T) \text{ for almost all }
 \boldsymbol{x} \in \Omega \}.
\end{align*}
For ease of notation, all inner products and norms in $L^2$ are denoted by
$(\cdot, \cdot)$ and $\| \cdot \|$,
if they are related to the whole space-time domain $Q_T$.
If they are associated with the spatial domain $\Omega$, then we write
$(\cdot, \cdot)_{\Omega}$
and $\| \cdot \|_{\Omega}$, which denote
the standard inner products and norms of the space $L^2(\Omega)$.
The symbols $(\cdot, \cdot)_{1,\Omega}$
and $\| \cdot \|_{1,\Omega}$ denote
the standard inner products and norms of $H^1(\Omega)$. 

The functions used in our analysis will typically be 
presented as Fourier series, i.e.,
\begin{align}
\label{def:FourierAnsatz}
 v(\boldsymbol{x},t) = v_0^c(\boldsymbol{x}) +
      \sum_{k=1}^{\infty} \left(v_k^c(\boldsymbol{x}) \cos(k \omega t)
                          + v_k^s(\boldsymbol{x}) \sin(k \omega t)\right)
\end{align}
with the Fourier coefficients 
\begin{align*}
 \begin{aligned}
 v_0^c(\boldsymbol{x}) &= \frac{1}{T} \int_0^T v(\boldsymbol{x},t) \,dt, \\
 v_k^c(\boldsymbol{x}) = \frac{2}{T} \int_0^T v(\boldsymbol{x},t) \cos(k \omega t)\,&dt, \qquad \qquad
 v_k^s(\boldsymbol{x}) = \frac{2}{T} \int_0^T v(\boldsymbol{x},t) \sin(k \omega t)\,dt,  
 \end{aligned}
\end{align*}
where $T$ and $\omega = 2 \pi /T$ denote the periodicity and the frequency, respectively.
Moreover, we define additional function spaces, see \cite{LRW:LangerWolfmayr:2013},
in order to derive a symmetric variational formulation of problem
(\ref{problem:parPeriodicBVP:PDE})-(\ref{problem:parPeriodicBVP:PerCond}).
The function spaces $H^{0,\frac{1}{2}}_{per}(Q_T)$, 
$H^{1,\frac{1}{2}}_{per}(Q_T)$ and $H^{1,\frac{1}{2}}_{0,per}(Q_T)$ are defined by
\begin{align*}
 H^{0,\frac{1}{2}}_{per}(Q_T) &= \{ u \in L^2(Q_T) : \big\| \partial^{1/2}_t u \big\|  < \infty \}, \\
 H^{1,\frac{1}{2}}_{per}(Q_T) &= \{ u \in H^{1,0}(Q_T) : \big\| \partial^{1/2}_t u \big\|  < \infty \}, \\
 H^{1,\frac{1}{2}}_{0,per}(Q_T) &= \{ u \in H^{1,\frac{1}{2}}_{per}(Q_T): u = 0 \mbox{ on } \Sigma_T \},
\end{align*}
respectively, where $\big\| \partial^{1/2}_t u \big\| $ is defined
in the Fourier space by the relation
\begin{align}
\label{definition:H01/2seminorm}
 \big\| \partial^{1/2}_t u \big\| ^2 :=
 |u|_{H^{0,\frac{1}{2}}(Q_T)}^2:= 
 \frac{T}{2} \sum_{k=1}^{\infty} k \omega \|\boldsymbol{u}_k\|_{\Omega}^2,
\end{align}
where $\boldsymbol{u}_k := (u_k^c, u_k^s)$ for all $k \in \mathbb{N}$.
These spaces are equipped with the scalar products
\begin{align}
\label{definition:H01/2products}
\begin{aligned}
 \big(\partial^{1/2}_t u, \partial^{1/2}_t v\big) :=
 \frac{T}{2} \sum_{k=1}^{\infty} k \omega (\boldsymbol{u}_k,\boldsymbol{v}_k)_{\Omega}, \qquad
 \big(\sigma
 \partial^{1/2}_t u, \partial^{1/2}_t v\big) :=
 \frac{T}{2} \sum_{k=1}^{\infty} k \omega (\sigma
 \boldsymbol{u}_k,\boldsymbol{v}_k)_{\Omega}.
\end{aligned}
\end{align}
The seminorm and the norm of the space $H^{1,\frac{1}{2}}_{per}(Q_T)$ are defined
by the relations
\begin{align*}
 |u|_{H^{1,\frac{1}{2}}(Q_T)}^2
 = \|\nabla u\| ^2 + \|\partial_t^{1/2} u\| ^2
 = T \, \|\nabla u_0^c\|_{\Omega}^2
 + \frac{T}{2} \sum_{k=1}^{\infty} \left(k\omega \|\boldsymbol{u}_k\|_{\Omega}^2
 + \|\nabla \boldsymbol{u}_k\|_{\Omega}^2\right)
\end{align*}
and 
\begin{align*}
 \|u\|_{H^{1,\frac{1}{2}}(Q_T)}^2 &= \|u\| ^2 + |u|_{H^{1,\frac{1}{2}}(Q_T)}^2 \\
 &= T \, (\|u_0^c\|_{\Omega}^2 + \|\nabla u_0^c\|_{\Omega}^2)
 + \frac{T}{2} \sum_{k=1}^{\infty} \left((1+k\omega) \|\boldsymbol{u}_k\|_{\Omega}^2
 + \|\nabla \boldsymbol{u}_k\|_{\Omega}^2\right),
\end{align*}
respectively.
Furthermore, we define
\begin{align}
 \label{definition:vperp}
 \begin{aligned}
 v^{\perp}(\boldsymbol{x},t) &:= 
 \sum_{k=1}^{\infty} 
 \left(- v_k^c(\boldsymbol{x}) \sin(k \omega t)
 + v_k^s(\boldsymbol{x}) \cos(k \omega t)\right) \\
 &= 
 \sum_{k=1}^{\infty}
 \underbrace{(v_k^s(\boldsymbol{x}),-v_k^c(\boldsymbol{x}))}_{=:(-\boldsymbol{v}_k^{\perp})^T} \cdot
                        \left( \begin{array}{l}
                          \cos(k \omega t) \\
                          \sin(k \omega t)
                        \end{array} \right).
 \end{aligned}
\end{align}
Note that the relation $\|\boldsymbol{u}_k^\perp\|^2_{\Omega}
 = \|\boldsymbol{u}_k\|^2_{\Omega}$ is valid.
\begin{lemma}
 \label{lemma:H11/2IdentiesAndOrthogonalities}
 The identities
 \begin{align}
 \label{equation:H11/2identities}
 \begin{aligned}
  \big(\sigma \partial_t^{1/2} u,\partial_t^{1/2} v \big)  =
  \big(\sigma \partial_t u,v^{\perp} \big)  \quad \mbox{ and } \quad 
  \big(\sigma \partial_t^{1/2} u,\partial_t^{1/2} v^{\perp} \big)  =
  \big(\sigma \partial_t u,v \big) 
 \end{aligned}
 \end{align}
 are valid 
 for all $u \in H^{0,1}_{per}(Q_T)$ and $v \in H^{0,\frac{1}{2}}_{per}(Q_T)$.
\begin{proof}
Using the definition of the $\sigma$-weighted scalar product in (\ref{definition:H01/2products})
and inserting the Fourier expansions of
\begin{align*}
 \partial_t u(\boldsymbol{x},t) :=
    \sum_{k=1}^{\infty} [k \omega \, u_k^s(\boldsymbol{x}) \cos(k \omega t)
    - k \omega \, u_k^c(\boldsymbol{x}) \sin(k \omega t)]
\end{align*}
as well as (\ref{definition:vperp})
into the inner products, we obtain
\begin{align*}
 \big(\sigma \partial_t^{1/2} u,\partial_t^{1/2} v \big)
 &= \frac{T}{2} \sum_{k=1}^{\infty} k \omega (\sigma \boldsymbol{u}_k,\boldsymbol{v}_k)_\Omega
 = \frac{T}{2} \sum_{k=1}^{\infty} k \omega (\sigma \boldsymbol{u}_k^\perp,\boldsymbol{v}_k^\perp)_\Omega \\
 &= \frac{T}{2} \sum_{k=1}^{\infty} k \omega (\sigma (-\boldsymbol{u}_k^\perp),(-\boldsymbol{v}_k^\perp))_\Omega
 = \big(\sigma \partial_t u,v^{\perp} \big)
\end{align*}
with $\boldsymbol{u}_k^{\perp} = (-u_k^s,u_k^c)^T$ for all $k \in \mathbb{N}$,
and
\begin{align*} 
 \big(\sigma \partial_t^{1/2} u,\partial_t^{1/2} v^{\perp} \big)
 = \frac{T}{2} \sum_{k=1}^{\infty} k \omega (\sigma \boldsymbol{u}_k,\boldsymbol{v}_k^\perp)_\Omega 
 = \frac{T}{2} \sum_{k=1}^{\infty} k \omega (\sigma (-\boldsymbol{u}_k^\perp),\boldsymbol{v}_k)_\Omega
 = \big(\sigma \partial_t u,v \big).
\end{align*}
\end{proof}
\end{lemma}
Hence, 
the following orthogonality relations hold:
\begin{align}
\label{equation:orthorelation}
\begin{aligned}
 &\big(\sigma \partial_t u,u\big)  = 0 \quad \mbox{ and } \quad
 (\sigma u^{\perp},u)  = 0 \qquad \forall \, u \in H^{0,1}_{per}(Q_T), \\
 &\big(\sigma \partial^{1/2}_t u,\partial^{1/2}_t u^{\perp}\big)  = 0
 \quad \mbox{ and } \quad \big(\nu \nabla u, \nabla u^{\perp}\big)  = 0
 \qquad \forall \, u \in H^{1,\frac{1}{2}}_{per}(Q_T),
\end{aligned}
\end{align}
where, e.g.,
\begin{align*}
 \big(\nu \nabla u, \nabla u^{\perp}\big)
 = \sum_{k=1}^{\infty}
 (\nu \nabla \boldsymbol{u}_k, \nabla \boldsymbol{u}_k^\perp)_\Omega
 = 0 \qquad \forall \, u \in H^{1,\frac{1}{2}}_{per}(Q_T)
\end{align*}
with
$\nabla \boldsymbol{u}_k := ((\nabla u_k^c)^T, (\nabla u_k^s)^T)^T$ and
$\nabla \boldsymbol{u}_k^\perp := (- (\nabla u_k^s)^T, (\nabla u_k^c)^T)^T$ for all $k \in \mathbb{N}$.
The identity
\begin{align}
 \label{identity:H01/2}
 \int_0^T \xi \, \partial_t^{1/2} v^\perp \, dt = - \int_0^T \partial_t^{1/2} \xi^\perp \, v \, dt
 \qquad \forall \, \xi, v \in H^{0,\frac{1}{2}}_{per}(Q_T)
\end{align}
is also defined in the Fourier space yielding the definitions
\begin{align}
\label{definition:L2productinH01/2}
 \big(\xi,\partial_t^{1/2} v\big) 
 := \frac{T}{2} \sum_{k=1}^{\infty} (k \omega)^{1/2} (\boldsymbol{\xi}_k,\boldsymbol{v}_k)_{\Omega}
\end{align}
as well as
\begin{align*}
 \partial_t^{1/2} \xi(\boldsymbol{x},t) &:=
 \sum_{k=1}^\infty (k \omega)^{1/2} \left(\xi_k^c(\boldsymbol{x})\cos(k \omega t)
 + \xi_k^s(\boldsymbol{x})\sin(k \omega t)\right)
\end{align*}
and
\begin{align*}
 \partial_t^{1/2} \xi^\perp(\boldsymbol{x},t) &:=
 \sum_{k=1}^\infty (k \omega)^{1/2} \left(-\xi_k^s(\boldsymbol{x})\cos(k \omega t)
 + \xi_k^c(\boldsymbol{x})\sin(k \omega t)\right).
\end{align*}
Hence,
\begin{align*}
 \big(\xi,\partial_t^{1/2} v^\perp\big) 
 &= \frac{T}{2} \sum_{k=1}^{\infty} (k \omega)^{1/2} (\boldsymbol{\xi}_k,\boldsymbol{v}_k^\perp)_{\Omega}
 = -\big(\partial_t^{1/2} \xi,v^\perp\big), \\
 \big(\xi,\partial_t^{1/2} v^\perp\big) 
 &= \frac{T}{2} \sum_{k=1}^{\infty} (k \omega)^{1/2} (\boldsymbol{\xi}_k,\boldsymbol{v}_k^\perp)_{\Omega}
 = \frac{T}{2} \sum_{k=1}^{\infty} (k \omega)^{1/2} (-\boldsymbol{\xi}_k^\perp,\boldsymbol{v}_k)_{\Omega} \\
 &= - \frac{T}{2} \sum_{k=1}^{\infty} (k \omega)^{1/2} (\boldsymbol{\xi}_k^\perp,\boldsymbol{v}_k)_{\Omega}
 = -\big(\partial_t^{1/2} \xi^\perp, v\big) 
\end{align*}
and all these
identities coincide with the
identities (\ref{equation:H11/2identities}) in Lemma~\ref{lemma:H11/2IdentiesAndOrthogonalities}.

We note that for functions presented in terms of Fourier series the standard
Friedrichs inequality holds in the form
\begin{align}
\label{inequality:Friedrichs:FourierSpace}
 \begin{aligned}
 \|\nabla u\| ^2 &= \int_{Q_T} |\nabla u|^2 \, d\boldsymbol{x}\,dt
 = T \, \|\nabla u_0^c\|_{\Omega}^2 + \frac{T}{2} \sum_{k=1}^\infty \|\nabla \boldsymbol{u}_k\|_{\Omega}^2 \\
 &\geq \frac{1}{C_F^2} \left(T \, \|u_0^c\|_{\Omega}^2 + \frac{T}{2} \sum_{k=1}^\infty \|\boldsymbol{u}_k\|_{\Omega}^2 \right) 
 = \frac{1}{C_F^2} \|u\| ^2.
 \end{aligned}
\end{align}

In order to derive the space-time variational formulation of the parabolic time-periodic problem
(\ref{problem:parPeriodicBVP:PDE})-(\ref{problem:parPeriodicBVP:PerCond}),
the parabolic partial differential equation (\ref{problem:parPeriodicBVP:PDE})
is multiplied
by a test function
$v \in H^{1,\frac{1}{2}}_{0,per}(Q_T)$,
integrated over the space-time cylinder $Q_T$, and after integration
by parts with respect to the space and time variables,
the following ``symmetric'' space-time variational formulation of the parabolic time-periodic 
boundary value problem
(\ref{problem:parPeriodicBVP:PDE})-(\ref{problem:parPeriodicBVP:PerCond})
is obtained:
Given $f \in L^2(Q_T)$,
find $u \in H^{1,\frac{1}{2}}_{0,per}(Q_T)$ such that
\begin{align}
\label{problem:STVFAPost}
 \begin{aligned}
 a(u,v) = \int_{Q_T} f(\boldsymbol{x},t)\,v(\boldsymbol{x},t) \, d\boldsymbol{x}\,dt
 \qquad \forall \, v \in H^{1,\frac{1}{2}}_{0,per}(Q_T)
 \end{aligned}
\end{align}
with the space-time bilinear form
\begin{align}
\label{definition:STBF}
 \begin{aligned}
  a(u,v) = \int_{Q_T} \Big( \sigma(\boldsymbol{x}) \partial_t^{1/2} u(\boldsymbol{x},t)
    \, \partial_t^{1/2} v^{\perp}(\boldsymbol{x},t)
  + \nu(\boldsymbol{x}) \nabla u(\boldsymbol{x},t) \cdot \nabla v(\boldsymbol{x},t) \Big) d\boldsymbol{x}
  \,dt, 
 \end{aligned}
\end{align}
where all functions are given in their Fourier series expansion in time, i.e.,
everything has to be understood in the sense of 
(\ref{definition:H01/2seminorm}) and (\ref{definition:H01/2products}).
In particular, this Fourier series
approach makes sense due to the time-periodicity condition (for $u$ and $v$).

\section{Multiharmonic finite element approximation}
\label{Sec3:MhFEApprox}

Inserting the Fourier series ansatz (\ref{def:FourierAnsatz}) into (\ref{problem:STVFAPost}) and
exploiting the orthogonality of the functions $\cos(k \omega t)$ and $\sin(k \omega t)$ 
with respect to the inner product $(\cdot,\cdot)_{L^2(0,T)}$,
we arrive at the following variational formulation corresponding to every single mode
$k \in \mathbb{N}$:
Given $\boldsymbol{f}_k \in (L^2(\Omega))^2$,
find $\boldsymbol{u}_k \in \mathbb{V} := V \times V = (H^1_0(\Omega))^2$ such that
\begin{align}
\label{problem:STVFk}
\begin{aligned}
 \int_{\Omega} \left(\nu(\boldsymbol{x}) \nabla \boldsymbol{u}_k(\boldsymbol{x})
    \cdot \nabla \boldsymbol{v}_k(\boldsymbol{x}) +
    k \omega \, \sigma(\boldsymbol{x}) \boldsymbol{u}_k(\boldsymbol{x})
    \cdot \boldsymbol{v}_k^{\perp}(\boldsymbol{x})\right) d\boldsymbol{x}
 = \int_{\Omega} \boldsymbol{f}_k(\boldsymbol{x}) \cdot \boldsymbol{v}_k(\boldsymbol{x}) \, d\boldsymbol{x}
\end{aligned}
\end{align}
for all $\boldsymbol{v}_k \in \mathbb{V}$.
In the case $k=0$, we obtain the following
variational formulation: Given $f_0^c \in L^2(\Omega)$, find $u_0^c \in V = H^1_0(\Omega)$
such that
\begin{align}
\label{problem:STVF0}
 \int_{\Omega} \nu(\boldsymbol{x}) \nabla u_0^c(\boldsymbol{x}) \cdot \nabla v_0^c(\boldsymbol{x})\,d\boldsymbol{x}
 = \int_{\Omega} f_0^c(\boldsymbol{x}) \, v_0^c(\boldsymbol{x})\,d\boldsymbol{x}
\end{align}
for all $v_0^c \in V$.
The variational problems (\ref{problem:STVFk}) and (\ref{problem:STVF0}) have a unique solution
due to the Babu\v{s}ka-Aziz theorem, see \cite{LRW:Wolfmayr:2014}.
In order to  solve these problems numerically,
the Fourier series are truncated at a finite index $N$ and
the unknown Fourier coefficients
$
 \boldsymbol{u}_k = (u_k^c, u_k^s)^T \in \mathbb{V}
$
are approximated
by finite element functions
$
 \boldsymbol{u}_{kh} = (u_{kh}^c, u_{kh}^s)^T \in \mathbb{V}_h = V_h \times V_h \subset \mathbb{V}.
$
Here,
$ 
 V_h = \text{span} \{\varphi_1, \dots, \varphi_n\}
$
with the standard nodal basis 
$\{\varphi_i(\boldsymbol{x}) = \varphi_{ih}(\boldsymbol{x}): i=1,2,\dots,n_h \}$,
and $h$ denotes the usual discretization parameter 
such that $n = n_h = \text{dim} V_h = O(h^{-d})$.
We use continuous, piecewise linear functions
on the finite elements on a regular triangulation $\mathcal{T}_h$
to construct the finite element subspace $V_h$ and its basis, see, e.g.,
\cite{LRW:Braess:2005, LRW:Ciarlet:1978, LRW:JungLanger:2013, LRW:Steinbach:2008}.
Under the assumptions (\ref{assumptions:sigmaNu:sigmaStrictlyPositive}),
we then obtain the following saddle point system
\begin{align}
 \label{problem:MhFEDiscretizedSTVFk}
 \left( \begin{array}{cc}
     k \omega M_{h,\sigma} & -K_{h,\nu} \\
     -K_{h,\nu} & -k \omega M_{h,\sigma} \end{array} \right) \left( \begin{array}{c}
     \underline{u}_k^s \\
     \underline{u}_k^c \end{array} \right) = \left( \begin{array}{c}
     -\underline{f}^c_k \\
     -\underline{f}^s_k \end{array} \right),
\end{align}
which has to be solved with respect to the nodal parameter vectors
$\underline{u}_k^s  = (u_{k,i}^s)_{i=1,\dots,n} \in \mathbb{R}^n $ 
and 
$\underline{u}_k^c  = (u_{k,i}^c)_{i=1,\dots,n} \in \mathbb{R}^n $
of the finite element approximations
\begin{align*}
u_{kh}^s(\boldsymbol{x}) = \sum_{i=1}^n u_{k,i}^s \, \varphi_i(\boldsymbol{x})
\quad \mbox{and} \quad
u_{kh}^c(\boldsymbol{x}) = \sum_{i=1}^n u_{k,i}^c \, \varphi_i(\boldsymbol{x})
\end{align*}
to the unknown Fourier coefficients
$u_k^s(\boldsymbol{x})$ and $u_k^c(\boldsymbol{x})$, respectively.
The matrices $K_{h,\nu}$ and $M_{h,\sigma}$ correspond to the weighted stiffness matrix and
weighted mass matrix, respectively. Their entries are computed by the formulas
\begin{align*}
\begin{aligned}
 K_{h,\nu}^{ij} = \int_{\Omega} \nu \, \nabla \varphi_i \cdot \nabla \varphi_j \,d\boldsymbol{x} \qquad \text{and} \qquad
 M_{h,\sigma}^{ij} = \int_{\Omega} \sigma \, \varphi_i \, \varphi_j \,d\boldsymbol{x}
\end{aligned}
\end{align*}
with $i,j = 1,\dots,n$, whereas
\begin{align*}
\begin{aligned}
 \underline{f}^c_k = \Big\lbrack \int_{\Omega} f^c_k \, \varphi_j \,d\boldsymbol{x} \Big\rbrack_{j=1,\dots,n}
 \quad \text{and} \quad
 \underline{f}^s_k = \Big\lbrack \int_{\Omega} f^s_k \, \varphi_j \,d\boldsymbol{x} \Big\rbrack_{j=1,\dots,n}.
\end{aligned}
\end{align*}
In the case $k=0$,
the following linear system
arising from the variational problem (\ref{problem:STVF0}) is obtained:
\begin{align}
 \label{problem:MhFEDiscretizedSTVF0}
     K_{h,\nu} \, \underline{u}_0^c = \underline{f}^c_0.
\end{align}
Fast and robust solvers for the linear systems
(\ref{problem:MhFEDiscretizedSTVFk}) and (\ref{problem:MhFEDiscretizedSTVF0})
can be found in 
\cite{LRW:KollmannKolmbauerLangerWolfmayrZulehner:2013,LRW:KrausWolfmayr:2013,
LRW:LangerWolfmayr:2013,LRW:Wolfmayr:2014}.
We use these solvers in order to obtain
the multiharmonic finite element approximation
\begin{align}
 \label{definition:MultiharmonicFEApproximationU}
 u_{N h}(\boldsymbol{x},t) = u_{0h}^c(\boldsymbol{x})
   + \sum_{k=1}^N \left(u_{kh}^c(\boldsymbol{x}) \cos(k \omega t)
                      + u_{kh}^s(\boldsymbol{x}) \sin(k \omega t)\right)
\end{align}
of the exact solution $u(\boldsymbol{x},t)$.
The next section is devoted to computable a posteriori estimates of the difference
between $u_{N h}$ and $u$.

\section{Functional a posteriori error estimates}
\label{Sec4:FunctionalAPostErrorEstimates}

First, we present 
inf-sup and sup-sup conditions
for the bilinear form (\ref{definition:STBF}).
\begin{lemma}
\label{lemma:STBFinfsupsupsup}
 The space-time bilinear form $a(\cdot,\cdot)$ defined by (\ref{definition:STBF})
 satisfies the following inf-sup and sup-sup conditions:
 \begin{align}
 \label{inequality:STBFinfsupsupsup}
  \mu_1 \|u\|_{H^{1,\frac{1}{2}}(Q_T)} \leq
  \sup_{0 \not= v \in H^{1,\frac{1}{2}}_{0,per}(Q_T)}
    \frac{a(u,v)}{\|v\|_{H^{1,\frac{1}{2}}(Q_T)}} \leq
  \mu_2 \|u\|_{H^{1,\frac{1}{2}}(Q_T)}
 \end{align}
 for all $u \in H^{1,\frac{1}{2}}_{0,per}(Q_T)$
 with positive constants
 $\mu_1 = \frac{1}{\sqrt{2}}\min\{\frac{\underline{\nu}}{C_F^2+1},\underline{\sigma}\}$ and
 $\mu_2 = \max\{\overline{\sigma},\overline{\nu}\}$,
 where $C_F$ is the constant coming from the Friedrichs inequality.
\begin{proof}
Using the triangle and Cauchy-Schwarz inequalities,
we obtain the estimate
 \begin{align*}
   |a(u,v)| &=
   \Big|\int_{Q_T} \left(\sigma(\boldsymbol{x}) \partial^{1/2}_t u(\boldsymbol{x},t) \,\partial^{1/2}_t v^{\perp}(\boldsymbol{x},t) 
  + \nu(\boldsymbol{x}) \nabla u(\boldsymbol{x},t)  \cdot \nabla v(\boldsymbol{x},t)  \right)\, d\boldsymbol{x}\,dt\Big| \\
   &\leq \overline{\sigma} \, \big\|\partial^{1/2}_t u\big\| \big\|\partial^{1/2}_t v\big\|
   + \overline{\nu} \, \|\nabla u\| \|\nabla v\|
   \leq \max\{\overline{\sigma},\overline{\nu}\} \, |u|_{H^{1,\frac{1}{2}}(Q_T)} |v|_{H^{1,\frac{1}{2}}(Q_T)} \\
   &\leq \mu_2 \, \|u\|_{H^{1,\frac{1}{2}}(Q_T)} \|v\|_{H^{1,\frac{1}{2}}(Q_T)}
 \end{align*}
with the constant $\mu_2 = \max\{\overline{\sigma},\overline{\nu}\}$,
which justifies the right hand-side inequality in (\ref{inequality:STBFinfsupsupsup}).

In order to
prove the left-hand side inequality, we select the test function
$v = u - u^{\perp}$
and estimate the supremum from below.
Using the $\sigma$- and $\nu$-weighted orthogonality relations
(\ref{equation:orthorelation}) and the
Friedrichs inequality
(\ref{inequality:Friedrichs:FourierSpace}),
we find that
\begin{align*}
   a(u,u) &=
   \int_{Q_T} \left(\sigma(\boldsymbol{x}) \partial^{1/2}_t u(\boldsymbol{x},t) \,\partial^{1/2}_t u^{\perp}(\boldsymbol{x},t) 
  + \nu(\boldsymbol{x}) \nabla u(\boldsymbol{x},t)  \cdot \nabla u(\boldsymbol{x},t)  \right)\, d\boldsymbol{x}\,dt \\
   &= \int_{Q_T} \nu(\boldsymbol{x}) \nabla u(\boldsymbol{x},t) \cdot \nabla u(\boldsymbol{x},t) \, d\boldsymbol{x}\,dt 
   \geq \underline{\nu} \int_{Q_T} |\nabla u|^2 \, d\boldsymbol{x}\,dt
   \geq \frac{\underline{\nu}}{c_F^2+1} \| u\|_{H^{1,0}(Q_T)}^2
\end{align*}
  and
  \begin{align*}
   a(u,-u^{\perp}) &=
   \int_{Q_T} \left(\sigma(\boldsymbol{x}) \partial^{1/2}_t u(\boldsymbol{x},t) \,\partial^{1/2}_t u(\boldsymbol{x},t) 
  - \nu(\boldsymbol{x}) \nabla u(\boldsymbol{x},t)  \cdot \nabla u^{\perp}(\boldsymbol{x},t)  \right)\, d\boldsymbol{x}\,dt \\
   &= \int_{Q_T} \sigma(\boldsymbol{x}) \partial^{1/2}_t u(\boldsymbol{x},t) \,\partial^{1/2}_t u(\boldsymbol{x},t)  \, d\boldsymbol{x}\,dt
   \geq \underline{\sigma} \, \big\| \partial^{1/2}_t u \big\|^2.
  \end{align*}
Combining these estimates, we have
  \begin{align*}
  \sup_{0 \not= v \in H^{1,\frac{1}{2}}_{0,per}(Q_T)}
    \frac{a(u,v)}{\|v\|_{H^{1,\frac{1}{2}}(Q_T)}} &\geq \frac{a(u,u-u^\perp)}{\|u-u^\perp\|_{H^{1,\frac{1}{2}}(Q_T)}}
   \geq \frac{\frac{\underline{\nu}}{c_F^2+1} \| u\|_{H^{1,0}(Q_T)}^2
   + \underline{\sigma} \, \big\| \partial^{1/2}_t u \big\|^2}{\|v\|_{H^{1,\frac{1}{2}}(Q_T)}} \\
   &\geq \frac{\min\{\frac{\underline{\nu}}{c_F^2+1},\underline{\sigma}\}
   \|u\|_{H^{1,\frac{1}{2}}(Q_T)}^2}{\sqrt{2} \|u\|_{H^{1,\frac{1}{2}}(Q_T)}}
   = \mu_1 \, \|u\|_{H^{1,\frac{1}{2}}(Q_T)},
  \end{align*}
  with the constant $\mu_1 = \frac{1}{\sqrt{2}}\min\{\frac{\underline{\nu}}{c_F^2+1},\underline{\sigma}\}$.
\end{proof}
\end{lemma}
\begin{remark}
 Since the condition $u=0$ is imposed on the whole boundary, 
 we can easily find an upper bound of $C_F$. Indeed, $C_F(\Omega) \leq C_F(\hat \Omega)$
 if $\hat \Omega \supset \Omega$.
 Since for such domains as rectangles or balls the Friedrichs constants are known,
 we can easily obtain an upper bound of
 $C_F$ for any Lipschitz domain.
\end{remark}
\begin{corollary}
 Since the norm $|\cdot|_{H^{1,\frac{1}{2}}(Q_T)}$ is equivalent to the norm
 $\|\cdot\|_{H^{1,\frac{1}{2}}(Q_T)}$ due to the Friedrichs inequality, the estimate
 (\ref{inequality:STBFinfsupsupsup}) implies
 \begin{align}
 \label{inequality:STBFinfsupsupsup:Seminorm}
  \tilde \mu_1 |u|_{H^{1,\frac{1}{2}}(Q_T)} \leq
  \sup_{0 \not= v \in H^{1,\frac{1}{2}}_{0,per}(Q_T)}
    \frac{a(u,v)}{|v|_{H^{1,\frac{1}{2}}(Q_T)}} \leq
  \tilde \mu_2 |u|_{H^{1,\frac{1}{2}}(Q_T)}
 \end{align}
 for all $u \in H^{1,\frac{1}{2}}_{0,per}(Q_T)$
 with positive constants $\tilde \mu_1 = \frac{1}{\sqrt{2}}\min\{\underline{\nu},\underline{\sigma}\}$ and
 $\tilde \mu_2 = \mu_2 = \max\{\overline{\sigma},\overline{\nu}\}$. 
\end{corollary}
We now move on to
the main part of this section related to a posteriori error estimation.
Let a function $\eta$ be an approximation of $u$.
First, we assume that $\eta$ is
a bit more regular than $u$.
More precisely, we set
$\eta \in H^{1,1}_{0,per}(Q_T)$.
This is of course true for the multiharmonic finite element approximation $u_{Nh}$,
which will later play the role of $\eta$.
Now, the ultimative goal 
is to deduce a computable upper bound of the error
$
 e := u - \eta
$
in $H^{1,\frac{1}{2}}_{0,per}(Q_T)$.
First, we notice that (\ref{problem:STVFAPost})
implies the integral identity
\begin{align}
 \label{problem:STVFAPost:Error}
  \begin{aligned}
 \int_{Q_T} &\Big( \sigma (\boldsymbol{x})
 \partial_t^{1/2} (u - \eta) \, \partial_t^{1/2} v^{\perp}
  + \nu (\boldsymbol{x})
  \nabla (u - \eta) \cdot
  \nabla v \Big) d\boldsymbol{x}\,dt \\
 &= \int_{Q_T} \Big(f\,v - \sigma (\boldsymbol{x})
 \partial_t^{1/2} \eta \, \partial_t^{1/2} v^{\perp}
  - \nu (\boldsymbol{x})
  \nabla \eta \cdot \nabla v \Big) \, d\boldsymbol{x}\,dt,
 \end{aligned}
\end{align}
which
is valid for all $v \in H^{1,\frac{1}{2}}_{0,per}(Q_T)$.
Here, the linear functional
\begin{align*}
 \mathcal{F}_\eta(v) := \int_{Q_T} \Big(f\,v - \sigma (\boldsymbol{x})
 \partial_t^{1/2} \eta \, \partial_t^{1/2} v^{\perp}
  - \nu (\boldsymbol{x}) \nabla \eta \cdot \nabla v \Big) \, d\boldsymbol{x}\,dt.
\end{align*}
is defined on $v \in H^{1,\frac{1}{2}}_{0,per}(Q_T)$.
Now, identity (\ref{problem:STVFAPost:Error}) can be rewritten in the form
\begin{align}
 \label{problem:STVFAPost:Error:Definition}
a(e,v) = \mathcal{F}_\eta(v).
\end{align}
Hence,
getting an upper bound of the error is reduced to finding the quantities
\begin{align}
 \label{inequality:supRHS}
 \sup_{0 \not= v \in H^{1,\frac{1}{2}}_{0,per}(Q_T)}
    \frac{\mathcal{F}_\eta(v)}{\|v\|_{H^{1,\frac{1}{2}}(Q_T)}}
 \qquad \qquad \text{ or } \qquad \qquad
 \sup_{0 \not= v \in H^{1,\frac{1}{2}}_{0,per}(Q_T)}
    \frac{\mathcal{F}_\eta(v)}{|v|_{H^{1,\frac{1}{2}}(Q_T)}}.
\end{align}
In order to find them, we reconstruct the functional
$\mathcal{F}_\eta(v)$ using
the identity 
\begin{align}
\label{equation:identityEtaH11/2}
 \big(\sigma \partial_t^{1/2} \eta,\partial_t^{1/2} v^{\perp} \big)  =
 \big(\sigma \partial_t \eta,v \big) \qquad \forall \, \eta \in H^{1,1}_{0,per}(Q_T)
 \quad \forall \, v \in H^{1,\frac{1}{2}}_{0,per}(Q_T),
\end{align}
which follows from (\ref{equation:H11/2identities})
and the identity
\begin{align*}
 \int_\Omega \text{div} \, \boldsymbol{\tau} \, v \, d\boldsymbol{x}
 = - \int_\Omega \boldsymbol{\tau} \cdot \nabla v \, d\boldsymbol{x},
\end{align*}
which is valid for any $v \in H^1_0(\Omega)$ and any
\begin{align*}
\boldsymbol{\tau} \in H(\text{div}_{\boldsymbol{x}},Q_T) := \{\boldsymbol{\tau} \in [L^2(Q_T)]^d 
 : \text{div}_{\boldsymbol{x}} \, 
 \boldsymbol{\tau}(\cdot,t) \in L^2(\Omega) \text{ for a.e. } t \in (0,T)
 \}.
\end{align*}
For ease of notation, the index $\boldsymbol{x}$ in $\text{div}_{\boldsymbol{x}}$ will be henceforth omitted,
i.e., $\text{div} = \text{div}_{\boldsymbol{x}}$ denotes the generalized
spatial divergence.
Using the Cauchy-Schwarz inequality leads to
\begin{align}
\label{inequality:supRHS:CS}
\begin{aligned}
 \mathcal{F}_\eta(v) &= \int_{Q_T} \Big(f\,v - \sigma (\boldsymbol{x})
 \partial_t \eta \, v + \text{div} \, \boldsymbol{\tau} \, v + (\boldsymbol{\tau} - \nu (\boldsymbol{x}) \nabla \eta)
 \cdot \nabla v \Big) \, d\boldsymbol{x}\,dt \\
 &\leq \|\mathcal{R}_1(\eta,\boldsymbol{\tau})\|  \|v\| 
 + \|\mathcal{R}_2(\eta,\boldsymbol{\tau})\|  \|\nabla v\|,
\end{aligned}
\end{align}
where
\begin{align*}
 \mathcal{R}_1(\eta,\boldsymbol{\tau})
   := \sigma \partial_t \eta - \text{div} \, \boldsymbol{\tau} - f \qquad \text{ and } \qquad
 \mathcal{R}_2(\eta,\boldsymbol{\tau})
   := \boldsymbol{\tau} - \nu \nabla \eta.
\end{align*}
In view of (\ref{inequality:Friedrichs:FourierSpace}),
we have
\begin{align*}
 \mathcal{F}_\eta(v) &\leq \|\mathcal{R}_1(\eta,\boldsymbol{\tau})\|  \|v\| 
 + \|\mathcal{R}_2(\eta,\boldsymbol{\tau})\|  \|\nabla v\| \\
 &\leq \|\mathcal{R}_1(\eta,\boldsymbol{\tau})\|  C_F \, \|\nabla v\| 
 + \|\mathcal{R}_2(\eta,\boldsymbol{\tau})\|  \|\nabla v\|  
 = \left(C_F \, \|\mathcal{R}_1(\eta,\boldsymbol{\tau})\| 
 + \|\mathcal{R}_2(\eta,\boldsymbol{\tau})\| \right) \|\nabla v\|.
\end{align*}
Hence, we obtain
\begin{align}
\label{inequality:aposteriorEstimateH11/2Seminorm0Step}
\begin{aligned}
  \sup_{0 \not= v \in H^{1,\frac{1}{2}}_{0,per}(Q_T)}
    \frac{\mathcal{F}_\eta(v)}{|v|_{H^{1,\frac{1}{2}}(Q_T)}}
  &\leq \sup_{0 \not= v \in H^{1,\frac{1}{2}}_{0,per}(Q_T)}
  \frac{\left(C_F \, \|\mathcal{R}_1(\eta,\boldsymbol{\tau})\| 
  + \|\mathcal{R}_2(\eta,\boldsymbol{\tau})\| \right)
    \|\nabla v\| }{|v|_{H^{1,\frac{1}{2}}(Q_T)}} \\
  &= \sup_{0 \not= v \in H^{1,\frac{1}{2}}_{0,per}(Q_T)}
  \frac{\left(C_F \, \|\mathcal{R}_1(\eta,\boldsymbol{\tau})\| 
  + \|\mathcal{R}_2(\eta,\boldsymbol{\tau})\| \right)
    \|\nabla v\| }{(\|\nabla v\| ^2 + \|\partial_t^{1/2} v\| ^2)^{1/2}} \\
  &\leq C_F \, \|\mathcal{R}_1(\eta,\boldsymbol{\tau})\| 
    + \|\mathcal{R}_2(\eta,\boldsymbol{\tau})\|.
\end{aligned}
\end{align}
We use (\ref{inequality:STBFinfsupsupsup:Seminorm}), i.e.,
  \begin{align*}
  |u - \eta|_{H^{1,\frac{1}{2}}(Q_T)} \leq \frac{1}{\tilde \mu_1}
  \sup_{0 \not= v \in H^{1,\frac{1}{2}}_{0,per}(Q_T)}
    \frac{a(u - \eta,v)}{|v|_{H^{1,\frac{1}{2}}(Q_T)}}
  = \frac{1}{\tilde \mu_1}
  \sup_{0 \not= v \in H^{1,\frac{1}{2}}_{0,per}(Q_T)}
    \frac{\mathcal{F}_\eta(v)}{|v|_{H^{1,\frac{1}{2}}(Q_T)}},
 \end{align*}
and arrive at
the following result:
\begin{theorem}
\label{theorem:aposteriorEstimateH11/2Seminorm}
 Let $\eta \in H^{1,1}_{0,per}(Q_T)$ and the bilinear form $a(\cdot,\cdot)$ satisfy
 (\ref{inequality:STBFinfsupsupsup:Seminorm}). Then,
 \begin{align}
  \label{inequality:aposteriorEstimateH11/2Seminorm}
  |u - \eta|_{H^{1,\frac{1}{2}}(Q_T)} \leq \frac{1}{\tilde \mu_1}
  \left(C_F \, \|\mathcal{R}_1(\eta,\boldsymbol{\tau})\| 
    + \|\mathcal{R}_2(\eta,\boldsymbol{\tau})\| \right)
  =: \mathcal{M}_{|\cdot|}^{\oplus}(\eta,\boldsymbol{\tau}),
 \end{align}
 where $\tilde \mu_1 = \frac{1}{\sqrt{2}}\min\{\underline{\nu},\underline{\sigma}\}$
 and $\boldsymbol{\tau} \in H(\emph{div},Q_T)$.
\end{theorem}

We can also deduce an upper bound of the full $H^{1,\frac{1}{2}}$-norm.
Indeed,
\begin{align*}
 \mathcal{F}_\eta(v)
 &\leq \|\mathcal{R}_1(\eta,\boldsymbol{\tau})\|  \|v\| 
 + \|\mathcal{R}_2(\eta,\boldsymbol{\tau})\|  \|\nabla v\|  \\
 &\leq \left(\|\mathcal{R}_1(\eta,\boldsymbol{\tau})\| ^2
   + \|\mathcal{R}_2(\eta,\boldsymbol{\tau})\| ^2\right)^{1/2}
  \left(\|v\| ^2 + \|\nabla v\| ^2\right)^{1/2}.
\end{align*}
In view of (\ref{inequality:STBFinfsupsupsup}), we obtain
\begin{align*}
  \sup_{0 \not= v \in H^{1,\frac{1}{2}}_{0,per}(Q_T)}
    \frac{\mathcal{F}_\eta(v)}{\|v\|_{H^{1,\frac{1}{2}}(Q_T)}}
  &\leq \sup_{0 \not= v \in H^{1,\frac{1}{2}}_{0,per}(Q_T)}
  \frac{\left(\|\mathcal{R}_1(\eta,\boldsymbol{\tau})\| ^2
    + \|\mathcal{R}_2(\eta,\boldsymbol{\tau})\| ^2\right)^{1/2}
  \left(\|v\| ^2 + \|\nabla v\| ^2\right)^{1/2}}{\|v\|_{H^{1,\frac{1}{2}}(Q_T)}} \\
  &\leq \left(\|\mathcal{R}_1(\eta,\boldsymbol{\tau})\| ^2
    + \|\mathcal{R}_2(\eta,\boldsymbol{\tau})\| ^2\right)^{1/2}.
\end{align*}
Altogether, we deduce a similar estimate for $\|e\|_{H^{1,\frac{1}{2}}(Q_T)}$.
\begin{theorem}
\label{theorem:aposteriorEstimateH11/2Norm}
 Let $\eta \in H^{1,1}_{0,per}(Q_T)$ and the bilinear form $a(\cdot,\cdot)$ satisfy
 (\ref{inequality:STBFinfsupsupsup}). Then,
 \begin{align}
 \label{inequality:aposteriorEstimateH11/2Norm}
  \|u - \eta\|_{H^{1,\frac{1}{2}}(Q_T)} \leq \frac{1}{\mu_1}
  \left(\|\mathcal{R}_1(\eta,\boldsymbol{\tau})\| ^2
    + \|\mathcal{R}_2(\eta,\boldsymbol{\tau})\| ^2\right)^{1/2}
  =: \mathcal{M}^\oplus_{\|\cdot\|}(\eta,\boldsymbol{\tau}),
 \end{align}
 where $\boldsymbol{\tau} \in H(\emph{div},Q_T)$ and now
 $\mu_1 = \frac{1}{\sqrt{2}}\min\{\frac{\underline{\nu}}{C_F^2+1},\underline{\sigma}\}$.
\end{theorem}
The functionals
$\mathcal{M}^\oplus_{|\cdot|}(\eta,\boldsymbol{\tau})$ and
$\mathcal{M}^\oplus_{\|\cdot\|}(\eta,\boldsymbol{\tau})$
present guaranteed and computable upper bounds (majorants) of the error
with respect to the
$H^{1,\frac{1}{2}}$-norm.
\begin{remark}
 It is easy to see that the majorants are nonnegative functionals vanishing if and only if
 $\eta = u$ and $\boldsymbol{\tau} = \nu \nabla u$.
 Indeed, if 
 $\mathcal{R}_1(\eta,\boldsymbol{\tau}) = 0$ and $\mathcal{R}_2(\eta,\boldsymbol{\tau}) = 0$, then
$
  \sigma \partial_t \eta - \emph{div} \, \boldsymbol{\tau} = f
$
and
$
  \boldsymbol{\tau} = \nu \nabla \eta.
$
 Since $\eta \in H^{1,1}_{0,per}(Q_T)$ is a periodic function and
 satisfies the Dirichlet condition on $\Sigma_T$, it is the solution.
 On the other hand,
 $\mathcal{R}_i(u,\nu \nabla u) = 0$, $i=1,2$.
\end{remark}

\subsection*{The multiharmonic approximation}

Since $f \in L^2(Q_T)$, it can be expanded into a Fourier series.
Moreover, we choose our approximation $\eta$ of the solution $u$
as well as the vector-valued function $\boldsymbol{\tau}$ to be
truncated Fourier series, i.e.,
\begin{align}
\label{definition:MhApproxEtaTau}
\begin{aligned}
 \eta(\boldsymbol{x},t) &= \eta_0^c(\boldsymbol{x}) + \sum_{k=1}^N \left(\eta_k^c(\boldsymbol{x}) \cos(k \omega t)
 + \eta_k^s(\boldsymbol{x}) \sin(k \omega t)\right), \\
 \boldsymbol{\tau}(\boldsymbol{x},t) &= \boldsymbol{\tau}_0^c(\boldsymbol{x}) + \sum_{k=1}^N \left(\boldsymbol{\tau}_k^c(\boldsymbol{x}) \cos(k \omega t)
 + \boldsymbol{\tau}_k^s(\boldsymbol{x}) \sin(k \omega t)\right),
\end{aligned}
\end{align}
where all Fourier coefficients 
are from the space $L^2(\Omega)$ and are defined by the relations
\begin{align*}
\begin{aligned}
 \eta_0^c(\boldsymbol{x}) &= \frac{1}{T} \int_0^T \eta(\boldsymbol{x},t) \, dt, \qquad \qquad \qquad 
 \boldsymbol{\tau}_0^c(\boldsymbol{x}) = \frac{1}{T} \int_0^T \boldsymbol{\tau}(\boldsymbol{x},t) \, dt, \\
 \eta_k^c(\boldsymbol{x}) &= \frac{2}{T} \int_0^T \eta(\boldsymbol{x},t) \cos(k \omega t) \, dt, \qquad
 \boldsymbol{\tau}_k^c(\boldsymbol{x}) = \frac{2}{T} \int_0^T \boldsymbol{\tau}(\boldsymbol{x},t) \cos(k \omega t) \, dt,\\
 \eta_k^s(\boldsymbol{x}) &= \frac{2}{T} \int_0^T \eta(\boldsymbol{x},t) \sin(k \omega t) \, dt, \qquad
 \boldsymbol{\tau}_k^s(\boldsymbol{x}) = \frac{2}{T} \int_0^T \boldsymbol{\tau}(\boldsymbol{x},t) \sin(k \omega t) \, dt.
\end{aligned}
\end{align*}
Hence, we get
\begin{align*}
 \partial_t \eta(\boldsymbol{x},t) &= \sum_{k=1}^N \left(k \omega \, \eta_k^s(\boldsymbol{x}) \cos(k \omega t)
    - k \omega \, \eta_k^c(\boldsymbol{x}) \sin(k \omega t)\right), \\
 \nabla \eta(\boldsymbol{x},t) &= \nabla \eta_0^c(\boldsymbol{x})
 + \sum_{k=1}^N \left(\nabla \eta_k^c(\boldsymbol{x}) \, \cos(k \omega t)
 + \nabla \eta_k^s(\boldsymbol{x}) \, \sin(k \omega t)\right), \\
 \text{div} \, \boldsymbol{\tau}(\boldsymbol{x},t) &= \text{div} \, \boldsymbol{\tau}_0^c(\boldsymbol{x})
 + \sum_{k=1}^N \left(\text{div} \, \boldsymbol{\tau}_k^c(\boldsymbol{x}) \, \cos(k \omega t)
 + \text{div} \, \boldsymbol{\tau}_k^s(\boldsymbol{x}) \, \sin(k \omega t)\right),
\end{align*}
and
the $L^2(Q_T)$-norms of the functions
\begin{align*}
 \mathcal{R}_1(\eta,\boldsymbol{\tau})
 = \sigma \partial_t \eta - \text{div} \, \boldsymbol{\tau} - f \qquad \text{ and } \qquad
 \mathcal{R}_2(\eta,\boldsymbol{\tau}) = \boldsymbol{\tau} - \nu \nabla \eta
\end{align*}
can easily be  computed.
Thus, we arrive at
\begin{align*}
 \|\mathcal{R}_1(\eta,\boldsymbol{\tau})\|^2
  = &\,T \|\text{div} \, \boldsymbol{\tau}_0^c+ f_0^c\|_{\Omega}^2 
      + \frac{T}{2} \sum_{k=1}^N
      \left(\|- k \omega \, \sigma \eta_k^s + \text{div} \, \boldsymbol{\tau}_k^c + f_k^c\|_{\Omega}^2
      + \|k \omega \, \sigma \eta_k^c + \text{div} \, \boldsymbol{\tau}_k^s + f_k^s\|_{\Omega}^2\right) \\
      &+ \frac{T}{2} \sum_{k=N+1}^\infty \left(\|f_k^c\|_{\Omega}^2 + \|f_k^s\|_{\Omega}^2\right) \\
  = &\,T \|\text{div} \, \boldsymbol{\tau}_0^c+ f_0^c\|_{\Omega}^2 + \frac{T}{2} \sum_{k=1}^N
      \|k \omega \, \sigma \boldsymbol{\eta}_k^\perp
      + \text{\textbf{div}} \, \boldsymbol{\tau}_k + \boldsymbol{f}_k\|_{\Omega}^2 
      + \frac{T}{2} \sum_{k=N+1}^\infty \|\boldsymbol{f}_k\|_{\Omega}^2,
\end{align*}
where $\boldsymbol{\eta}_k^\perp = (-\eta_k^s,\eta_k^c)^T$,
$\text{\textbf{div}} \, \boldsymbol{\tau}_k
= (\text{div} \, \boldsymbol{\tau}_k^c,\text{div} \, \boldsymbol{\tau}_k^s)^T$,
and
\begin{align*}
 \|\mathcal{R}_2(\eta,\boldsymbol{\tau})\| ^2
 &= \int_{Q_T} \left|\boldsymbol{\tau} - \nu \nabla \eta\right|^2 d\boldsymbol{x}\,dt \\
 &= T \|\boldsymbol{\tau}_0^c - \nu \nabla \eta_0^c\|_{\Omega}^2
    + \frac{T}{2} \sum_{k=1}^N
      \left(\|\boldsymbol{\tau}_k^c - \nu \nabla \eta_k^c\|_{\Omega}^2
      + \|\boldsymbol{\tau}_k^s - \nu \nabla \eta_k^s\|_{\Omega}^2\right) \\
 &= T \|\boldsymbol{\tau}_0^c - \nu \nabla \eta_0^c\|_{\Omega}^2 + \frac{T}{2} \sum_{k=1}^N
      \|\boldsymbol{\tau}_k - \nu \nabla \boldsymbol{\eta}_k\|_{\Omega}^2
\end{align*}
with $\boldsymbol{\tau}_k = ((\boldsymbol{\tau}_k^c)^T,(\boldsymbol{\tau}_k^s)^T)^T$.
\begin{remark}
 We note that the remainder term
\begin{align*}
 \mathcal{E}_N := \frac{T}{2} \sum_{k=N+1}^\infty \|\boldsymbol{f}_k\|_{\Omega}^2
 = \frac{T}{2} \sum_{k=N+1}^\infty \left(\|f_k^c\|_{\Omega}^2 + \|f_k^s\|_{\Omega}^2\right)
\end{align*} 
 is always computable, due to the knowledge on the given data $f$.
 In some cases, the computation of $\mathcal{E}_N$ is very
 easy, for example,
 if $f$ is multiharmonic.
 However, even in the most complicated cases, in which $f = f(\boldsymbol{x},t)$ and
 we do not refer to special (e.g., extra regularity) properties, the term
 $\mathcal{E}_N$ can be precomputed as $\|f - f_N\|$,
 where $f_N$ is the truncated Fourier series of $f$.
\end{remark}
In fact, the $L^2$-norms of $\mathcal{R}_1$ and $\mathcal{R}_2$ corresponding to
every single mode $k$ are decoupled.
Altogether, it follows that
\begin{align*}
 \begin{aligned}
 \|\mathcal{R}_1(\eta,\boldsymbol{\tau})\|^2
 = &\,T \|{\mathcal{R}_1}^c_0(\boldsymbol{\tau}_0^c)\|_{\Omega}^2
      + \frac{T}{2} \sum_{k=1}^N
      \left(\|{\mathcal{R}_1}^c_k(\eta_k^s,\boldsymbol{\tau}_k^c)\|_{\Omega}^2
      + \|{\mathcal{R}_1}^s_k(\eta_k^c,\boldsymbol{\tau}_k^s)\|_{\Omega}^2\right) + \mathcal{E}_N, \\
  \|\mathcal{R}_2(\eta,\tau)\|^2
  = &\,T \|{\mathcal{R}_2}^c_0(\eta_0^c,\boldsymbol{\tau}_0^c)\|_{\Omega}^2 
    + \frac{T}{2} \sum_{k=1}^N
      \left(\|{\mathcal{R}_2}^c_k(\eta_k^c,\boldsymbol{\tau}_k^c)\|_{\Omega}^2
      + \|{\mathcal{R}_2}^s_k(\eta_k^s,\boldsymbol{\tau}_k^s)\|_{\Omega}^2\right),
 \end{aligned}
\end{align*}
where  
${\mathcal{R}_1}^c_0(\boldsymbol{\tau}_0^c) := \text{div} \, \boldsymbol{\tau}_0^c + f_0^c$,\;
${\mathcal{R}_2}^c_0(\eta_0^c,\boldsymbol{\tau}_0^c) := \boldsymbol{\tau}_0^c - \nu \nabla \eta_0^c$,
and, for $k=1,\dots,N$, we have
\begin{align}
\label{definition:R1R2:Fourier}
 \begin{aligned}
 {\mathcal{R}_1}^c_k(\eta_k^s,\boldsymbol{\tau}_k^c) &:= - k \omega \, \sigma \eta_k^s
 + \text{div} \, \boldsymbol{\tau}_k^c + f_k^c,\\
 {\mathcal{R}_1}^s_k(\eta_k^c,\boldsymbol{\tau}_k^s) &:= k \omega \, \sigma \eta_k^c
 + \text{div} \, \boldsymbol{\tau}_k^s + f_k^s,\\
 {\mathcal{R}_2}^c_k(\eta_k^c,\boldsymbol{\tau}_k^c) &:= \boldsymbol{\tau}_k^c - \nu \nabla \eta_k^c,\\
{\mathcal{R}_2}^s_k(\eta_k^s,\boldsymbol{\tau}_k^s) &:= \boldsymbol{\tau}_k^s - \nu \nabla \eta_k^s.
 \end{aligned}
\end{align}
\begin{corollary}
The error majorants $\mathcal{M}_{|\cdot|}^{\oplus}(\eta,\boldsymbol{\tau})$ and
$\mathcal{M}_{\|\cdot\|}^{\oplus}(\eta,\boldsymbol{\tau})$ can be presented in the forms
\begin{align*}
  \mathcal{M}_{|\cdot|}^{\oplus}(\eta,\boldsymbol{\tau})
  = \frac{1}{\tilde \mu_1}
  &\Big(C_F \, \|\mathcal{R}_1(\eta,\boldsymbol{\tau})\| 
    + \|\mathcal{R}_2(\eta,\boldsymbol{\tau})\| \Big) \\
  = \frac{1}{\tilde \mu_1}
  &\Big(C_F \, \big(T \|{\mathcal{R}_1}^c_0(\boldsymbol{\tau}_0^c)\|_{\Omega}^2
      + \frac{T}{2} \sum_{k=1}^N
      \left(\|{\mathcal{R}_1}^c_k(\eta_k^s,\boldsymbol{\tau}_k^c)\|_{\Omega}^2
      + \|{\mathcal{R}_1}^s_k(\eta_k^c,\boldsymbol{\tau}_k^s)\|_{\Omega}^2\right) + \mathcal{E}_N
      \big)^{1/2} \\
      &+ \Big(T \|{\mathcal{R}_2}^c_0(\eta_0^c,\boldsymbol{\tau}_0^c)\|_{\Omega}^2
    +  \frac{T}{2} \sum_{k=1}^N
      \left(\|{\mathcal{R}_2}^c_k(\eta_k^c,\boldsymbol{\tau}_k^c)\|_{\Omega}^2
      + \|{\mathcal{R}_2}^s_k(\eta_k^s,\boldsymbol{\tau}_k^s)\|_{\Omega}^2\right)\Big)^{1/2} \Big)
\end{align*}
and 
\begin{align*}
 \mathcal{M}_{\|\cdot\|}^{\oplus}(\eta,\boldsymbol{\tau}) = \frac{1}{\mu_1}
  &\Big(\|\mathcal{R}_1(\eta,\boldsymbol{\tau})\| ^2
  + \|\mathcal{R}_2(\eta,\boldsymbol{\tau})\| ^2\Big)^{1/2} \\
  = \frac{1}{\mu_1}
  &\Big(T \big(\|{\mathcal{R}_1}^c_0(\boldsymbol{\tau}_0^c)\|_{\Omega}^2
    + \|{\mathcal{R}_2}^c_0(\eta_0^c,\boldsymbol{\tau}_0^c)\|_{\Omega}^2\big)
      + \frac{T}{2} \sum_{k=1}^N
      \big(\|{\mathcal{R}_1}^c_k(\eta_k^s,\boldsymbol{\tau}_k^c)\|_{\Omega}^2 \\
      &+ \|{\mathcal{R}_1}^s_k(\eta_k^c,\boldsymbol{\tau}_k^s)\|_{\Omega}^2
      + \|{\mathcal{R}_2}^c_k(\eta_k^c,\boldsymbol{\tau}_k^c)\|_{\Omega}^2
      + \|{\mathcal{R}_2}^s_k(\eta_k^s,\boldsymbol{\tau}_k^s)\|_{\Omega}^2\big)
    + \mathcal{E}_N
    \Big)^{1/2},
\end{align*}
where $\tilde \mu_1 = \frac{1}{\sqrt{2}}\min\{\underline{\nu},\underline{\sigma}\}$
and $\mu_1 = \frac{1}{\sqrt{2}} \min\{\frac{\underline{\nu}}{C_F^2+1},\underline{\sigma}\}$.
\end{corollary}
\begin{remark}
 Since the error (with respect to the truncation index $N$) between the exact solution $u$
 and its multiharmonic approximation $\eta$ decreases with $\mathcal{O}(N^{-1})$,
 see \cite{LRW:LangerWolfmayr:2013, LRW:Wolfmayr:2014}, the contributions in the majorants
 coming from the functionals ${\mathcal{R}_1}_k^c$ and ${\mathcal{R}_1}_k^s$ cannot blow up.
\end{remark}
We see that the majorants consist of computable quantities related to each
harmonic. Therefore, they not only evaluate the overall error, but also provide
an information on errors associated with a certain harmonic.
Moreover, since the respective quantities are integrals
over $\Omega$, their integrands serve
as indicators of spatial errors. Thus, the majorants contain a rich amount
of information to be utilized in various
adaptive procedures.
\begin{remark}
Let $f$ has a multiharmonic representation, i.e.,
 \begin{align*}
  f(\boldsymbol{x},t) = f_0^c(\boldsymbol{x}) + \sum_{k=1}^{N_f} \left(f_k^c(\boldsymbol{x}) \cos(k \omega t)
  + f_k^s(\boldsymbol{x}) \sin(k \omega t)\right),
 \end{align*}
where $N_f \in \mathbb{N}$ is defined by $f$.
If $N \geq N_f$, then 
$\eta$ is the exact solution of problem (\ref{problem:STVFAPost}) and
$\boldsymbol{\tau}$ is the exact flux if and only if the error majorants 
vanish, i.e.,
\begin{align}
 \label{equation:R12FourierVanish}
  \begin{aligned}
   {\mathcal{R}_1}_k^c = 0 \qquad \text{and} \qquad
   {\mathcal{R}_2}_k^c = 0 \qquad \forall \, k=0,1,\dots,N_f, \\
   {\mathcal{R}_1}_k^s = 0 \qquad \text{and} \qquad
   {\mathcal{R}_2}_k^s = 0 \qquad \forall \, k=1,2,\dots,N_f.
  \end{aligned}
\end{align}
Indeed, let the error majorants vanish.
Then, we deduce that
$- \emph{div} \, \boldsymbol{\tau}_0^c = f_0^c$ and $\boldsymbol{\tau}_0^c = \nu \nabla \eta_0^c$,
and furthermore we have 
$k \omega \, \sigma \eta_k^s - \emph{div} \, \boldsymbol{\tau}_k^c = f_k^c$,
$- k \omega \, \sigma \eta_k^c - \emph{div} \, \boldsymbol{\tau}_k^s = f_k^s$,
$\boldsymbol{\tau}_k^c = \nu \nabla \eta_k^c$ 
and  
$\boldsymbol{\tau}_k^s = \nu \nabla \eta_k^s$
for all $k=1,\dots,N_f$.
Therefore, collecting the harmonics, we find that
 \begin{align*}
  \boldsymbol{\tau}(\boldsymbol{x},t) &= \boldsymbol{\tau}_0^c(\boldsymbol{x})
  + \sum_{k=1}^{N_f} \left(\boldsymbol{\tau}_k^c(\boldsymbol{x}) \cos(k \omega t)
  + \boldsymbol{\tau}_k^s(\boldsymbol{x}) \sin(k \omega t)\right), \\
 \eta(\boldsymbol{x},t) &= \eta_0^c(\boldsymbol{x})
  + \sum_{k=1}^{N_f} \left(\eta_k^c(\boldsymbol{x}) \cos(k \omega t)
  + \eta_k^s(\boldsymbol{x}) \sin(k \omega t)\right)
 \end{align*}
and
 \begin{align*}
  \sigma \partial_t \eta - \emph{div} \, \boldsymbol{\tau} = f, \qquad
  \boldsymbol{\tau} = \nu \nabla \eta.
 \end{align*}
Since $\eta$ satisfies the boundary conditions and the equation, we conclude that $\eta = u$.
\end{remark}
Another approach to derive a majorant is to insert the Fourier series ansatz directly
to the bilinear form $a(u - \eta,v)$
and into the functional $\mathcal{F}_\eta(v)$ as defined 
in (\ref{problem:STVFAPost:Error}).
Then, we obtain the following integral identities associated with every mode:
\begin{align}
\label{problem:BFak:Error}
\begin{aligned}
  \int_{\Omega}
  &\left(\nu(\boldsymbol{x}) \nabla (\boldsymbol{u}_k(\boldsymbol{x})-\boldsymbol{\eta}_k(\boldsymbol{x}))
  \cdot \nabla \boldsymbol{v}_k(\boldsymbol{x}) +
  k \omega \, \sigma(\boldsymbol{x}) (\boldsymbol{u}_k(\boldsymbol{x})-\boldsymbol{\eta}_k(\boldsymbol{x}))
  \cdot \boldsymbol{v}_k^{\perp}(\boldsymbol{x})\right) d\boldsymbol{x} \\
  &= \int_{\Omega}
  \left(\boldsymbol{f}_k(\boldsymbol{x}) \cdot \boldsymbol{v}_k(\boldsymbol{x})
  - \nu(\boldsymbol{x}) \nabla \boldsymbol{\eta}_k(\boldsymbol{x})
  \cdot \nabla \boldsymbol{v}_k(\boldsymbol{x})  -
  k \omega \, \sigma(\boldsymbol{x}) \boldsymbol{\eta}_k(\boldsymbol{x})
  \cdot \boldsymbol{v}_k^{\perp}(\boldsymbol{x})\right) d\boldsymbol{x},
\end{aligned}
\end{align}
which are valid for all $\boldsymbol{v}_k \in (H^1_0(\Omega))^2$.
In the case $k=0$, the integral identity 
\begin{align}
\label{problem:BFa0:Error}
\begin{aligned}
  \int_{\Omega} \nu(\boldsymbol{x})\nabla (u_0^c(\boldsymbol{x})-\eta_0^c(\boldsymbol{x}))
  \cdot \nabla v_0^c(\boldsymbol{x})\,d\boldsymbol{x}
  = \int_{\Omega} \left(f_0^c(\boldsymbol{x}) \, v_0^c(\boldsymbol{x})
  -\nu(\boldsymbol{x})\nabla \eta_0^c(\boldsymbol{x})
  \cdot \nabla v_0^c(\boldsymbol{x})\right)\,d\boldsymbol{x}
\end{aligned}
\end{align}
is valid for all $v_0^c \in H^1_0(\Omega)$.
We define the left hand sides of (\ref{problem:BFak:Error}) and (\ref{problem:BFa0:Error}) by
\begin{align*}
 a_k(\boldsymbol{u}_k-\boldsymbol{\eta}_k,\boldsymbol{v}_k) \qquad \text{and} \qquad
 a_0(u_0^c-\eta_0^c,v_0^c),
\end{align*}
and the right hand sides by
\begin{align*}
 \mathcal{F}_{\boldsymbol{\eta}_k}(\boldsymbol{v}_k) \qquad \text{and} \qquad
 \mathcal{F}_{\eta_0^c}(v_0^c),
\end{align*}
respectively.
Let us start with the case $k=1,\dots,N$. Hence, 
an upper bound for the errors
 $\boldsymbol{e}_k := \boldsymbol{u}_k - \boldsymbol{\eta}_k$
in $(H^1_0(\Omega))^2$
has to be computed.
The bilinear form $a_k(\cdot,\cdot)$
meets the inf-sup condition
\begin{align}
\label{inequality:infsup:Error:k}
  \sup_{0 \not= \boldsymbol{v}_k \in (H^1_0(\Omega))^2}
    \frac{a_k(\boldsymbol{u}_k - \boldsymbol{\eta}_k,\boldsymbol{v}_k)}{\|\boldsymbol{v}_k\|_{1,\Omega}}
    \geq \underline{c}^k
    \, \|\boldsymbol{u}_k - \boldsymbol{\eta}_k\|_{1,\Omega}
\end{align}
with the inf-sup constant
$\underline{c}^k = \min\{\underline{\nu}, k \omega \, \underline{\sigma}\}/\sqrt{2}$.
By the same method as before, we reform the error functionals and obtain estimates for
\begin{align*}
 \sup_{0 \not= \boldsymbol{v}_k \in (H^1_0(\Omega))^2}
    \frac{\mathcal{F}_{\boldsymbol{\eta}_k}(\boldsymbol{v}_k)}{\|\boldsymbol{v}_k\|_{1,\Omega}}.
\end{align*}    
We introduce a collection of vector-valued functions
\begin{align*}
 \boldsymbol{\tau}_k = (\boldsymbol{\tau}_k^c,\boldsymbol{\tau}_k^s)^T, \qquad
 \boldsymbol{\tau}_k^c, \boldsymbol{\tau}_k^s \in H(\text{div},\Omega) :=
 \{\boldsymbol{\tau} \in [L^2(\Omega)]^d
 : \text{div} \, \boldsymbol{\tau} \in L^2(\Omega)\}
\end{align*}
and use the integral relations
\begin{align*}
\int_\Omega \text{div} \, \boldsymbol{\tau} \, v \, d\boldsymbol{x}
 = - \int_\Omega \boldsymbol{\tau} \cdot \nabla v \, d\boldsymbol{x}
 \qquad \forall \, v \in H^1_0(\Omega).
\end{align*}
It is easy to see that
\begin{align}
\label{inequality:supRHS:CS:k}
\begin{aligned}
 \mathcal{F}_{\boldsymbol{\eta}_k}(\boldsymbol{v}_k)
 = &\int_{\Omega}
  \big(\boldsymbol{f}_k \cdot \boldsymbol{v}_k
  - k \omega \, \sigma(\boldsymbol{x}) \boldsymbol{\eta}_k
  \cdot \boldsymbol{v}_k^{\perp}
  + \text{\textbf{div}} \, \boldsymbol{\tau}_k \cdot \boldsymbol{v}_k \\
  &+ \boldsymbol{\tau}_k \cdot \nabla \boldsymbol{v}_k
  - \nu(\boldsymbol{x}) \nabla \boldsymbol{\eta}_k
  \cdot \nabla \boldsymbol{v}_k \big) d\boldsymbol{x} \\
 = &\int_{\Omega}
  \big((\boldsymbol{f}_k
  + k \omega \, \sigma(\boldsymbol{x}) \boldsymbol{\eta}_k^\perp
  + \text{\textbf{div}} \, \boldsymbol{\tau}_k ) \cdot \boldsymbol{v}_k
  + (\boldsymbol{\tau}_k - \nu(\boldsymbol{x}) \nabla \boldsymbol{\eta}_k)
  \cdot \nabla \boldsymbol{v}_k \big) d\boldsymbol{x} \\
 \leq &\,\|{\mathcal{R}_1}_k(\boldsymbol{\eta}_k,\boldsymbol{\tau}_k)\|_{\Omega}
   \|\boldsymbol{v}_k\|_{\Omega}
 + \|{\mathcal{R}_2}_k(\boldsymbol{\eta}_k,\boldsymbol{\tau}_k)\|_{\Omega}
   \|\nabla \boldsymbol{v}_k\|_{\Omega} \\
 \leq &\,\left(\|{\mathcal{R}_1}_k(\boldsymbol{\eta}_k,\boldsymbol{\tau}_k)\|_{\Omega}^2
 + \|{\mathcal{R}_2}_k(\boldsymbol{\eta}_k,\boldsymbol{\tau}_k)\|_{\Omega}^2\right)^{1/2}
 \|\boldsymbol{v}_k\|_{1,\Omega},
\end{aligned}
\end{align}
where
\begin{align*}
 {\mathcal{R}_1}_k(\boldsymbol{\eta}_k,\boldsymbol{\tau}_k)
  = k \omega \, \sigma \boldsymbol{\eta}_k^\perp
  + \text{\textbf{div}} \, \boldsymbol{\tau}_k + \boldsymbol{f}_k
  &= (-k \omega \, \sigma \eta_k^s + \text{div} \, \boldsymbol{\tau}_k^c + f_k^c,
     k \omega \, \sigma \eta_k^c + \text{div} \, \boldsymbol{\tau}_k^s + f_k^s)^T \\
  &= ({\mathcal{R}_1}^c_k(\eta_k^s,\boldsymbol{\tau}_k^c),{\mathcal{R}_1}^s_k(\eta_k^c,\boldsymbol{\tau}_k^s))^T
\end{align*}
and
\begin{align*}
 {\mathcal{R}_2}_k(\boldsymbol{\eta}_k,\boldsymbol{\tau}_k)
  = \boldsymbol{\tau}_k - \nu \nabla \boldsymbol{\eta}_k
  &= (\boldsymbol{\tau}_k^c - \nu \nabla \eta_k^c,
     \boldsymbol{\tau}_k^s - \nu \nabla \eta_k^s)^T 
  = ({\mathcal{R}_2}^c_k(\eta_k^c,\boldsymbol{\tau}_k^c),{\mathcal{R}_2}^s_k(\eta_k^s,\boldsymbol{\tau}_k^s))^T.
\end{align*}
Hence, we have derived
the same results as in (\ref{definition:R1R2:Fourier})
for every mode $k=1,\dots,N$.
Using the estimate (\ref{inequality:supRHS:CS:k})
together with the inf-sup condition (\ref{inequality:infsup:Error:k}), we finally
arrive at the following upper bounds for every single mode $k=1,\dots,N$:
\begin{theorem}
\label{corollary:Fourierk:aposterioriEst:Norm}
Let $\boldsymbol{\eta}_k \in (H^1_0(\Omega))^2$ and the bilinear form $a_k(\cdot,\cdot)$
satisfy (\ref{inequality:infsup:Error:k}).
Then,
\begin{align}
 \label{inequality:aposterioriEstimateH11/2Norm:Fourierk}
  \|\boldsymbol{u}_k - \boldsymbol{\eta}_k\|_{1,\Omega}
  \leq \frac{1}{\underline{c}^k}
  \left(\|{\mathcal{R}_1}_k(\boldsymbol{\eta}_k,\boldsymbol{\tau}_k)\|_{\Omega}^2
  + \|{\mathcal{R}_2}_k(\boldsymbol{\eta}_k,\boldsymbol{\tau}_k)\|_{\Omega}^2\right)^{1/2}
  =: \mathcal{M}^{\oplus_k}_{\|\cdot\|}(\boldsymbol{\eta}_k,\boldsymbol{\tau}_k),
\end{align}
where $\underline{c}^k
= \min\{\underline{\nu}, k \omega \, \underline{\sigma}\}/\sqrt{2}$
and $\boldsymbol{\tau}_k = (\boldsymbol{\tau}_k^c, \boldsymbol{\tau}_k^s)^T$
with
$\boldsymbol{\tau}_k^c, \boldsymbol{\tau}_k^s \in H(\text{\emph{div}},\Omega)$.
\end{theorem}
Using the inf-sup condition
\begin{align}
\label{inequality:infsup:Error:k:Seminorm}
 \begin{aligned}
  \sup_{0 \not= \boldsymbol{v}_k \in (H^1_0(\Omega))^2}
    \frac{a_k(\boldsymbol{u}_k,\boldsymbol{v}_k)}{|\boldsymbol{v}_k|_{1,\Omega}}
  &= \sup_{0 \not= \boldsymbol{v}_k \in (H^1_0(\Omega))^2}
    \frac{\left(\nu \nabla \boldsymbol{u}_k, \nabla \boldsymbol{v}_k\right)_{\Omega} +
  k \omega \left(\sigma \boldsymbol{u}_k,
  \boldsymbol{v}_k^{\perp}\right)_{\Omega}}{|\boldsymbol{v}_k|_{1,\Omega}} \\
  &\geq
    \frac{\left(\nu \nabla \boldsymbol{u}_k,
    \nabla (\boldsymbol{u}_k-\boldsymbol{u}_k^\perp)\right)_{\Omega} +
  k \omega \left(\sigma \boldsymbol{u}_k,
  (\boldsymbol{u}_k-\boldsymbol{u}_k^\perp)^{\perp}\right)_{\Omega}}{
  |\boldsymbol{u}_k-\boldsymbol{u}_k^\perp|_{1,\Omega}} \\
  &= \frac{\left(\nu \nabla \boldsymbol{u}_k,
    \nabla \boldsymbol{u}_k\right)_{\Omega} +
  k \omega \left(\sigma \boldsymbol{u}_k,\boldsymbol{u}_k\right)_{\Omega}}{
  \sqrt{2} \, |\boldsymbol{u}_k|_{1,\Omega}}
  \geq \frac{\underline{\nu} \|\nabla \boldsymbol{u}_k\|_{\Omega}^2 +
  k \omega \,  \underline{\sigma} \|\boldsymbol{u}_k\|_{\Omega}^2}{
  \sqrt{2} \, |\boldsymbol{u}_k|_{1,\Omega}} \\
  &\geq \frac{\min\{\underline{\nu}, 
  k \omega \,  \underline{\sigma}\} 
  \|\boldsymbol{u}_k\|_{1,\Omega}^2}{
  \sqrt{2} \, |\boldsymbol{u}_k|_{1,\Omega}}
    \geq
    \frac{\min\{\underline{\nu}, 
    k \omega \, \underline{\sigma}
    \}}{\sqrt{2}}
    \,
    |\boldsymbol{u}_k|_{1,\Omega}
 \end{aligned}   
\end{align}
together with the estimate
\begin{align*}
\begin{aligned}
 \mathcal{F}_{\boldsymbol{\eta}_k}(\boldsymbol{v}_k)
 &\leq \|{\mathcal{R}_1}_k(\boldsymbol{\eta}_k,\boldsymbol{\tau}_k)\|_{\Omega}
   \|\boldsymbol{v}_k\|_{\Omega}
 + \|{\mathcal{R}_2}_k(\boldsymbol{\eta}_k,\boldsymbol{\tau}_k)\|_{\Omega}
   \|\nabla \boldsymbol{v}_k\|_{\Omega} \\
 &\leq \left(C_F \|{\mathcal{R}_1}_k(\boldsymbol{\eta}_k,\boldsymbol{\tau}_k)\|_{\Omega}
 + \|{\mathcal{R}_2}_k(\boldsymbol{\eta}_k,\boldsymbol{\tau}_k)\|_{\Omega}\right)
   |\boldsymbol{v}_k|_{1,\Omega}
\end{aligned}
\end{align*}
yields the following error majorant for $|\cdot|_{1,\Omega}$
with the same inf-sup constant $\underline{c}^k$:
\begin{theorem}
\label{corollary:Fourierk:aposterioriEst:Seminorm}
Let $\boldsymbol{\eta}_k \in (H^1_0(\Omega))^2$ and the bilinear form $a_k(\cdot,\cdot)$
satisfy (\ref{inequality:infsup:Error:k:Seminorm}).
Then,
\begin{align}
 \label{inequality:aposterioriEstimateH11/2Seminorm:Fourierk}
 \begin{aligned}
  |\boldsymbol{u}_k - \boldsymbol{\eta}_k|_{1,\Omega}
  &\leq
  \frac{1}{\underline{c}^k
  }
  \left(C_F \|{\mathcal{R}_1}_k(\boldsymbol{\eta}_k,\boldsymbol{\tau}_k)\|_{\Omega}
 + \|{\mathcal{R}_2}_k(\boldsymbol{\eta}_k,\boldsymbol{\tau}_k)\|_{\Omega}\right) 
  =: \mathcal{M}^{\oplus_k}_{|\cdot|}(\boldsymbol{\eta}_k,\boldsymbol{\tau}_k),
 \end{aligned}
\end{align}
where
$\underline{c}^k
= \min\{\underline{\nu},k \omega \, \underline{\sigma}
\}/\sqrt{2}$
and
$\boldsymbol{\tau}_k = (\boldsymbol{\tau}_k^c, \boldsymbol{\tau}_k^s)^T$
with
$\boldsymbol{\tau}_k^c, \boldsymbol{\tau}_k^s \in H(\text{\emph{div}},\Omega)$.
\end{theorem}
Now, we consider the case $k=0$. Here, an upper bound for the error
 $e_0^c := u_0^c - \eta_0^c$
in $H^1_0(\Omega)$ has to be computed.
The inf-sup condition
\begin{align}
\label{inequality:infsup:Error:k0}
  \sup_{0 \not= v_0^c \in H^1_0(\Omega)}
    \frac{a_0(u_0^c - \eta_0^c,v_0^c)}{\|v_0^c\|_{1,\Omega}}
    \geq \underline{c}^0_{\|\cdot\|} \, \|u_0^c - \eta_0^c\|_{1,\Omega}
\end{align}
with the inf-sup constant
$\underline{c}^0_{\|\cdot\|} = \underline{\nu}/(C_F^2+1)$
can be proved quite analogously to (\ref{inequality:infsup:Error:k}).
Moreover, one can easily show that
\begin{align}
\label{inequality:infsup:Error:k0:Seminorm}
  \sup_{0 \not= v_0^c \in H^1_0(\Omega)}
    \frac{a_0(u_0^c - \eta_0^c,v_0^c)}{|v_0^c|_{1,\Omega}}
  \geq \frac{a_0(u_0^c - \eta_0^c,u_0^c - \eta_0^c)}{|u_0^c - \eta_0^c|_{1,\Omega}}
    \geq \underline{c}^0_{|\cdot|} \, |u_0^c - \eta_0^c|_{1,\Omega} 
\end{align}
with $\underline{c}^0_{|\cdot|} = \underline{\nu}$,
since $\nu$ satisfies the assumptions (\ref{assumptions:sigmaNu:sigmaStrictlyPositive}).
By arguments similar to those used above for the modes $k$, we deduce the following estimates:
\begin{align}
 \label{inequality:aposteriorEstimateH11/2Norm:Fourierk0}
  \|u_0^c - \eta_0^c\|_{1,\Omega}
  \leq \frac{1}{\underline{c}^0_{\|\cdot\|}}
  \left(\|{\mathcal{R}_1}_0^c(\boldsymbol{\tau}_0^c)\|_{\Omega}^2
  + \|{\mathcal{R}_2}_0^c(\eta_0^c,\boldsymbol{\tau}_0^c)\|_{\Omega}^2 \right)^{1/2}
  =: \mathcal{M}^{\oplus_0}_{\|\cdot\|}(\eta_0^c,\boldsymbol{\tau}_0^c) 
\end{align}
and
\begin{align}
 \label{inequality:aposteriorEstimateH11/2Norm:Fourierk0:Seminorm}
  |u_0^c - \eta_0^c|_{1,\Omega} 
  \leq \frac{1}{\underline{c}^0_{|\cdot|}}
  \left(C_F \|{\mathcal{R}_1}_0^c(\boldsymbol{\tau}_0^c)\|_{\Omega}
  + \|{\mathcal{R}_2}_0^c(\eta_0^c,\boldsymbol{\tau}_0^c)\|_{\Omega} \right)
  =: \mathcal{M}^{\oplus_0}_{|\cdot|}(\eta_0^c,\boldsymbol{\tau}_0^c),
\end{align}
where 
$\boldsymbol{\tau}_0^c \in H(\text{div},\Omega)$, 
${\mathcal{R}_1}_0^c(\boldsymbol{\tau}_0^c) = f_0^c + \text{div} \, \boldsymbol{\tau}_0^c$
and
${\mathcal{R}_2}_0^c(\eta_0^c,\boldsymbol{\tau}_0^c) = \boldsymbol{\tau}_0^c - \nu \nabla \eta_0^c$.

\section{Numerical results}
\label{Sec5:NumericalResults}

In this section, we present and discuss results of numerical experiments
on computing functional a posteriori error estimates in the context of parabolic
time-periodic boundary value problems discretized by the MhFEM.
First, we present a numerical example with a given time-harmonic source term.
In the second example, we consider a given time-periodic, but not time-harmonic source term.
The computational domain $\Omega = (0,1) \times (0,1)$ is uniformly
decomposed into triangles, and standard continuous, piecewise linear
finite elements are used for the discretization in space. In this case, the
Friedrichs constant is $C_F = 1/(\sqrt{2}\pi)$.
In these two numerical experiments, we choose $\sigma = \nu = 1$.

The construction of $\eta$ and $\boldsymbol{\tau}$ is an important issue 
in order to obtain sharp guaranteed bounds from the majorants
$\mathcal{M}^{\oplus}_{\|\cdot\|}$ or $\mathcal{M}^{\oplus}_{|\cdot|}$.
As it has been already discussed in Section~\ref{Sec4:FunctionalAPostErrorEstimates}, 
we can choose multiharmonic finite element approximations
(\ref{definition:MhApproxEtaTau}) for $\eta$ and $\boldsymbol{\tau}$.
However, since the Fourier coefficients of $\eta$ are constructed by continuous,
piecewise linear approximations,
their gradients are only piecewise constant.
Then, $\nabla \eta_k^c, \nabla \eta_k^s \in L^2(\Omega)$, but
$\nabla \eta_k^c, \nabla \eta_k^s \not\in H(\text{div},\Omega)$,
$k = 1,\dots,N$.
Hence, a flux reconstruction is needed in order to obtain a
suitable flux $\boldsymbol{\tau} \in H(\text{div},Q_T)$.
A good reconstruction of the flux is an important and nontrivial topic.
We can regularize $\boldsymbol{\tau}$ by a post-processing operator 
which maps the $L^2$-functions into $H(\text{div},Q_T)$,
see \cite{LRW:Repin:2008}.
There are various techniques for realizing these post-processing steps such as, e.g.,
local post-processing by an elementwise averaging procedure
or by using Raviart-Thomas elements, see 
\cite{LRW:Repin:2008, LRW:MaliNeittaanmaekiRepin:2014} and the references therein.
In our numerical experiments, we use Raviart-Thomas elements of the lowest order,
see, e.g., \cite{LRW:RaviartThomas:1977, LRW:BrezziFortin:1991, LRW:RobertsThomas:1991}.
First, we define the normal fluxes on interior edges $E_{mn}$ by
\begin{align*}
 (\boldsymbol{\tau}^c_k \cdot n_{E_{mn}})|_{E_{mn}} &= (\lambda_{mn} (\nabla \eta_k^c)|_{T_m}
 + (1-\lambda_{mn}) (\nabla \eta_k^c)|_{T_n})\cdot n_{E_{mn}},
\\
 (\boldsymbol{\tau}^s_k \cdot n_{E_{mn}})|_{E_{mn}} &= (\lambda_{mn} (\nabla \eta_k^s)|_{T_m}
 + (1-\lambda_{mn}) (\nabla \eta_k^s)|_{T_n})\cdot n_{E_{mn}},
\end{align*}
for all $k = 1,\dots,N$,
with $\lambda_{mn} = 1/2$ due to uniform discretization.
Here, $(\nabla \eta_k^c)|_{T_m}$, $(\nabla \eta_k^s)|_{T_m}$, $(\nabla \eta_k^c)|_{T_n}$ and $(\nabla \eta_k^s)|_{T_n}$
are constant vectors on two arbitrary, neighboring elements $T_m$ and $T_n$.
On boundary edges, the only one existing flux is used.
Hence, three normal fluxes are defined on the three sides of each element.
Inside, we reconstruct the fluxes $\boldsymbol{\tau}_k = (\boldsymbol{\tau}^c_k,\boldsymbol{\tau}^s_k)^T$
by the standard lowest-order Raviart-Thomas ($RT^0$-) extension of normal fluxes with
\begin{align*}
RT^0(\mathcal{T}_h) := \{&\boldsymbol{\tau} \in (L^2(T))^2:
\forall \, T \in \mathcal{T}_h \quad \exists \, a,b,c \in \mathbb{R}
\quad \forall \, \boldsymbol{x} \in T, \\
&\boldsymbol{\tau}(\boldsymbol{x}) = (a,b)^T + c \, \boldsymbol{x} \text{ and }
[\boldsymbol{\tau}]_E \cdot n_E = 0 \,\, \forall 
\text{ interior edges } E \},
\end{align*}
where $[\boldsymbol{\tau}]_E$ denotes the jump of $\boldsymbol{\tau}$ across the edge $E$ shared by two neighboring
elements on a triangulation $\mathcal{T}_h$.
Altogether, it follows an averaged flux from $H(\text{div},\Omega)$, i.e.,
\begin{align*}
 \boldsymbol{\tau}^c_k = G_{RT}(\nabla \eta^c_k), \quad
 \boldsymbol{\tau}^s_k = G_{RT}(\nabla \eta^s_k), \quad
 G_{RT}:L^2(\Omega) \rightarrow H(\text{div},\Omega).
\end{align*}

In order to solve the saddle point systems (\ref{problem:MhFEDiscretizedSTVFk}) for $k=1,\dots,N$,
we use the AMLI preconditioner proposed by Kraus and Wolfmayr in \cite{LRW:KrausWolfmayr:2013}
with a proper $3$-refinement of the mesh as presented in \cite{LRW:KrausWolfmayr:2013}
for an inexact realization of the block-diagonal preconditioner
\begin{align}
\label{definition:robustPreconditionerForSTVFk}
 \mathcal{P} =
 \left( \begin{array}{cc}
     k \omega M_{h,\sigma} + K_{h,\nu} & 0 \\
     0 & k \omega M_{h,\sigma} + K_{h,\nu} \end{array} \right)
\end{align}
in the MINRES method.
The preconditioner (\ref{definition:robustPreconditionerForSTVFk})
was presented and discussed in \cite{LRW:Wolfmayr:2014}.
Here, we want to emphasize that the AMLI preconditioned MINRES solver
is robust and of optimal complexity, see \cite{LRW:KrausWolfmayr:2013, LRW:Wolfmayr:2014}.
This can be also observed in the numerical results of this paper.
We mention that, in all tables where the number of MINRES iterations
$n^{iter}_{\text{\tiny{MINRES}}}$ or of AMLI iterations $n^{iter}_{\text{\tiny{AMLI}}}$
is presented,
the iteration was stopped after reducing the initial residual by a factor of $10^{-6}$.
In each MINRES iteration step, we have used the AMLI preconditioner
according to \cite{LRW:KrausWolfmayr:2013}
with $8$ inner iterations.
The presented CPU times in seconds $t^{\text{sec}}$ include the
computational times for computing the majorants, which are very small in comparison to the
computational times of the solver.
All computations were performed on a PC with Intel(R) Xeon(R) CPU W3680 @ 3.33GHz.

In the {\bf first example}, we consider a given time-harmonic source term
\begin{align*}
 f(\boldsymbol{x},t) = 2(x_1(1-x_1)+x_2(1-x_2)) \cos(t) + x_1(1-x_1)x_2(x_2-1) \sin(t),
\end{align*}
where $T=2 \pi / \omega$ with $\omega = 1$. Hence, the Fourier coefficients of $f$ are simply given by
\begin{align*}
 f^c(\boldsymbol{x}) = 2(x_1(1-x_1)+x_2(1-x_2)), \qquad
 f^s(\boldsymbol{x}) = x_1(1-x_1)x_2(x_2-1),
\end{align*}
and we have to consider only one single mode $k=1$.
For simplicity, we omit the index $k$ in Example~1.
The exact solution is given by
\begin{align*}
 u(\boldsymbol{x},t) = x_1(x_1-1)x_2(x_2-1) \cos(t).
\end{align*}
Table~\ref{tab:Ex1} presents the number of MINRES iterations $n^{iter}_{\text{\tiny{MINRES}}}$,
the CPU times in seconds $t^{\text{sec}}$, the norms of $\mathcal{R}_1$ and $\mathcal{R}_2$, i.e.,
\begin{align*}
\|\mathcal{R}_{1}\|_{\Omega}^2 &= \|\mathcal{R}_1^c(\eta^s,\boldsymbol{\tau}^c)\|_{\Omega}^2
      + \|\mathcal{R}_1^s(\eta^c,\boldsymbol{\tau}^s)\|_{\Omega}^2 \\
      &= \|- \eta^s + \text{div} \, \boldsymbol{\tau}^c + f^c\|_{\Omega}^2
      + \|\eta^c + \text{div} \, \boldsymbol{\tau}^s + f^s\|_{\Omega}^2, \\
\|\mathcal{R}_{2}\|_{\Omega}^2 &= \|\mathcal{R}_2^c(\eta^s,\boldsymbol{\tau}^c)\|_{\Omega}^2
      + \|\mathcal{R}_2^s(\eta^c,\boldsymbol{\tau}^s)\|_{\Omega}^2
      = \|\boldsymbol{\tau}^c - \nabla \eta^c\|_{\Omega}^2
      +  \|\boldsymbol{\tau}^s - \nabla \eta^s\|_{\Omega}^2,
\end{align*}
as well as the majorants
\begin{align*}
 \mathcal{M}^{\oplus}_{|\cdot|} 
 = \frac{1}{\tilde \mu_1 
 } \big(C_F \|\mathcal{R}_{1}\|_{\Omega} + \|\mathcal{R}_{2}\|_{\Omega} \big),
\end{align*}
where $\tilde \mu_1 = 1/\sqrt{2}$, and the corresponding efficiency indices
\begin{align}
 I_{\text{eff}} =
 \frac{\mathcal{M}^{\oplus}_{|\cdot|}}{|\boldsymbol{u}-\boldsymbol{\eta}|_{1,\Omega}},
\end{align}
obtained on grids of different mesh sizes. Here,
$\boldsymbol{u} = \boldsymbol{u}(\boldsymbol{x}) = (u^c(\boldsymbol{x}),u^s(\boldsymbol{x}))^T$ denotes the
vector of the exact solution's
Fourier coefficients $u^c(\boldsymbol{x}) = x_1(x_1-1)x_2(x_2-1)$ and $u^s(\boldsymbol{x}) = 0$.
\begin{table}[!ht]
\begin{center}
\begin{tabular}{|c|cccccc|}
  \hline
   grid & $n^{iter}_{\text{\tiny{MINRES}}}$ & $t^{\text{sec}}$ & $\|\mathcal{R}_{1}\|_{\Omega}$
   & $\|\mathcal{R}_{2}\|_{\Omega}$ & $\mathcal{M}^{\oplus}_{|\cdot|}$ & $I_{\text{eff}}$  \\
  \hline
   $9 \times 9$     & 14 & 0.00  & 1.657e-01 & 4.604e-03 & 5.926e-02 & 1.976 \\ 
   $27 \times 27$   & 14 & 0.03  & 6.381e-02 & 8.313e-05 & 2.043e-02 & 1.583 \\ 
   $81 \times 81$   & 12 & 0.24  & 2.186e-02 & 7.545e-06 & 6.968e-03 & 1.530 \\ 
   $243 \times 243$ & 12 & 2.43  & 7.334e-03 & 5.155e-07 & 2.335e-03 & 1.504 \\ 
   $729 \times 729$ & 12 & 22.25 & 2.449e-03 & 3.298e-08 & 7.797e-04 & 1.498 \\ 
  \hline
\end{tabular}
\end{center}
\caption{Majorant and its parts (Example 1).}
\label{tab:Ex1}
\end{table}

In Table~\ref{tab:Ex1}, we observe the robustness and optimality of the
AMLI preconditioned MINRES method as presented in \cite{LRW:KrausWolfmayr:2013, LRW:Wolfmayr:2014}.
More precisely, the computational times increase with a factor of nine that exactly reveals
the optimal computational complexity of the method according to the $3$-refinement of the mesh.
One can see that the norms of $\mathcal{R}_1$ reduce as a factor of three and
the norms of $\mathcal{R}_2$ even better than as a factor of nine.
Hence, the applied flux reconstruction is efficient.
Altogether, the majorant reduces as a factor of three by trisection of the mesh size and is 
of the same order of convergence as of the exact error measured in the $H^1(\Omega)$-seminorm.
This is also observed in the efficiency index that is already quite small on the
$27 \times 27$-mesh and decreases up to a value of $1.498$ on the (finest) $729 \times 729$-mesh.

In the {\bf second example}, we consider a given time-analytic, but not time-harmonic source term
\begin{align*}
 f(\boldsymbol{x},t) = e^t \sin^2(t) \sin(x_1 \pi) \sin(x_2 \pi)((1+2\pi^2)\sin(t)+3\cos(t)),
\end{align*}
where $T=2 \pi / \omega$ with $\omega = 1$.
The exact solution is given by
\begin{align*}
 u(\boldsymbol{x},t) = e^t \sin^3(t) \sin(x_1 \pi) \sin(x_2 \pi).
\end{align*}
The Fourier coefficients of the Fourier series expansion of the source term $f$ in time
can be computed analytically.
We truncate the Fourier series and approximate the Fourier coefficients by
finite element functions as it was presented before.
Then, we solve the systems (\ref{problem:MhFEDiscretizedSTVFk}) and (\ref{problem:MhFEDiscretizedSTVF0})
for all $k \in \{0,\dots,N\}$ with $N=8$, reconstruct the fluxes by a $RT^0$-extension and
then compute the corresponding majorants.
Table~\ref{tab:Ex2:0} presents the number of AMLI iterations $n^{iter}_{\text{\tiny{AMLI}}}$,
the CPU times in seconds $t^{\text{sec}}$, the norms of ${\mathcal{R}_1}_0^c$ and ${\mathcal{R}_2}_0^c$, i.e.,
\begin{align*}
\|{\mathcal{R}_1}_0^c\|_{\Omega}^2 &= \|\text{div} \, \boldsymbol{\tau}_0^c + f_0^c\|_{\Omega}^2, \qquad \qquad
\|{\mathcal{R}_2}_0^c\|_{\Omega}^2 = \|\boldsymbol{\tau}_0^c - \nabla \eta_0^c\|_{\Omega}^2,
\end{align*}
as well as the majorants $\mathcal{M}^{\oplus_0}_{|\cdot|}$ as presented in
(\ref{inequality:aposteriorEstimateH11/2Norm:Fourierk0:Seminorm})
with $c^0_{|\cdot|} = \underline{\nu} = 1$, and
the corresponding efficiency indices
\begin{align}
 I_{\text{eff}}^0 =
 \frac{\mathcal{M}^{\oplus_0}_{|\cdot|}}{|u_0^c - \eta_0^c|_{1,\Omega}}
\end{align}
obtained on grids of different mesh sizes.
\begin{table}[!ht]
\begin{center}
\begin{tabular}{|c|cccccc|}
  \hline
   grid & $n^{iter}_{\text{\tiny{AMLI}}}$ & $t^{\text{sec}}$ & $\|{\mathcal{R}_1}_0^c\|_{\Omega}$
   & $\|{\mathcal{R}_2}_0^c\|_{\Omega}$ & $\mathcal{M}^{\oplus_0}_{|\cdot|}$ & $I_{\text{eff}}^0$  \\
  \hline
   $9 \times 9$     & 21 & 0.00 & 6.317e+01 & 1.773e+00 & 1.599e+01 & 1.315 \\ 
   $27 \times 27$   & 23 & 0.00 & 2.349e+01 & 3.796e-02 & 5.325e+00 & 1.064 \\ 
   $81 \times 81$   & 23 & 0.03 & 7.927e+00 & 2.865e-03 & 1.787e+00 & 1.020 \\ 
   $243 \times 243$ & 22 & 0.27 & 2.646e+00 & 1.886e-04 & 5.957e-01 & 1.006 \\ 
   $729 \times 729$ & 22 & 2.45 & 8.821e-01 & 1.183e-05 & 1.986e-01 & 1.002 \\ 
  \hline
\end{tabular}
\end{center}
\caption{Majorant $\mathcal{M}^{\oplus_0}_{|\cdot|}$ and its parts
(Example 2).}
\label{tab:Ex2:0}
\end{table}

For $k=0$, one has to solve the system (\ref{problem:MhFEDiscretizedSTVF0}).
We observe in Table~\ref{tab:Ex2:0} that the AMLI solver
presented by Kraus and Wolfmayr in \cite{LRW:KrausWolfmayr:2013}
is of optimal computational complexity and the efficiency decreases up to a value of
$1.002$.
Moreover, Tables~\ref{tab:Ex2:1} -- \ref{tab:Ex2:8}
present the number of MINRES iterations $n^{iter}_{\text{\tiny{MINRES}}}$,
the CPU times in seconds $t^{\text{sec}}$, the norms of ${\mathcal{R}_1}_k$ and ${\mathcal{R}_2}_k$, i.e.,
\begin{align*}
\|{\mathcal{R}_1}_k\|_{\Omega}^2 &= \|{\mathcal{R}_1}_k^c(\eta_k^s,\boldsymbol{\tau}_k^c)\|_{\Omega}^2
      + \|{\mathcal{R}_1}_k^s(\eta_k^c,\boldsymbol{\tau}_k^s)\|_{\Omega}^2 \\
      &= \|- k \omega \, \eta_k^s + \text{div} \, \boldsymbol{\tau}_k^c + f_k^c\|_{\Omega}^2
      + \|k \omega \, \eta_k^c + \text{div} \, \boldsymbol{\tau}_k^s + f_k^s\|_{\Omega}^2, \\
\|{\mathcal{R}_2}_k\|_{\Omega}^2 &= \|{\mathcal{R}_2}_k^c(\eta_k^s,\boldsymbol{\tau}_k^c)\|_{\Omega}^2
      + \|{\mathcal{R}_2}_k^s(\eta_k^c,\boldsymbol{\tau}_k^s)\|_{\Omega}^2
      = \|\boldsymbol{\tau}_k^c - \nabla \eta_k^c\|_{\Omega}^2
      +  \|\boldsymbol{\tau}_k^s - \nabla \eta_k^s\|_{\Omega}^2,
\end{align*}
as well as the majorants
$\mathcal{M}^{\oplus_k}_{|\cdot|}$ as presented in
(\ref{inequality:aposterioriEstimateH11/2Seminorm:Fourierk}) 
with $\underline{c}^k = \min\{\underline{\nu},k \omega \, \underline{\sigma}\}/\sqrt{2} = 1/\sqrt{2}$
for $k \in \{1,\dots,8\}$,
and, finally, the corresponding efficiency indices
\begin{align}
 I_{\text{eff}}^k =
 \frac{\mathcal{M}^{\oplus_k}_{|\cdot|}}{|\boldsymbol{u}_k-\boldsymbol{\eta}_k|_{1,\Omega}}
\end{align}
obtained on grids of different mesh sizes.

The results of Tables~\ref{tab:Ex2:1} -- \ref{tab:Ex2:8} regarding the number of MINRES iterations
$n^{iter}_{\text{\tiny{MINRES}}}$ and the computational times are all similar and can be compared
to our Example~1. Moreover, the reduction factors of $\|{\mathcal{R}_1}_k\|_{\Omega}$,
$\|{\mathcal{R}_2}_k\|_{\Omega}$ and $\mathcal{M}^{\oplus_k}_{|\cdot|}$ as well as the values
of the efficiency indices $I_{\text{eff}}^k$ are approximately the same. This demonstrates
the robustness of the method with respect to the modes $k$ and the accurateness of
the majorants $\mathcal{M}^{\oplus_k}_{|\cdot|}$. Moreover, the values of $\|{\mathcal{R}_1}_k\|_{\Omega}$,
$\|{\mathcal{R}_2}_k\|_{\Omega}$ and $\mathcal{M}^{\oplus_k}_{|\cdot|}$ decrease for increasing $k$.
This is also illustrated in Table~\ref{tab:Ex2:global}.
In this table,
we finally compare the results from Tables~\ref{tab:Ex2:1} -- \ref{tab:Ex2:8}
that were computed on the $729 \times 729$-mesh.
Hence, the
results computed on the $729 \times 729$-mesh are again presented
for all $k \in \{0,\dots,8\}$, and, then, for the overall functional error estimates.
Here, the error majorant is given by
\begin{align*}
  \mathcal{M}_{|\cdot|}^{\oplus}(\eta,\boldsymbol{\tau})
  = \frac{1}{\tilde \mu_1}
  &\Big(C_F \, \|\mathcal{R}_1(\eta,\boldsymbol{\tau})\| 
    + \|\mathcal{R}_2(\eta,\boldsymbol{\tau})\| \Big) \\
  = \frac{1}{\tilde \mu_1}
  &\Big(C_F \, \big(T \|{\mathcal{R}_1}^c_0(\boldsymbol{\tau}_0^c)\|_{\Omega}^2
      + \frac{T}{2} \sum_{k=1}^N
      \left(\|{\mathcal{R}_1}^c_k(\eta_k^s,\boldsymbol{\tau}_k^c)\|_{\Omega}^2
      + \|{\mathcal{R}_1}^s_k(\eta_k^c,\boldsymbol{\tau}_k^s)\|_{\Omega}^2\right) + \mathcal{E}_N
      \big)^{1/2} \\
      &+ \Big(T \|{\mathcal{R}_2}^c_0(\eta_0^c,\boldsymbol{\tau}_0^c)\|_{\Omega}^2
    +  \frac{T}{2} \sum_{k=1}^N
      \left(\|{\mathcal{R}_2}^c_k(\eta_k^c,\boldsymbol{\tau}_k^c)\|_{\Omega}^2
      + \|{\mathcal{R}_2}^s_k(\eta_k^s,\boldsymbol{\tau}_k^s)\|_{\Omega}^2\right)\Big)^{1/2} \Big),
\end{align*}
where $\tilde \mu_1 = 1/\sqrt{2}$, and the remainder term 
\begin{align*}
 \mathcal{E}_N = \frac{T}{2} \sum_{k=N+1}^\infty \|\boldsymbol{f}_k\|_{\Omega}^2
 = \frac{T}{2} \sum_{k=N+1}^\infty \left(\|f_k^c\|_{\Omega}^2 + \|f_k^s\|_{\Omega}^2\right)
\end{align*}
has to be computed in order to get $\|\mathcal{R}_1\|$. Remember that the remainder term
can be precomputed exactly as $\|f-f_N\|$, since $f$ is the given data and $f_N$ its truncated
Fourier series.
Altogether, we obtain a global efficiency index of $1.404$ on the $729 \times 729$-mesh.
\begin{table}[!ht]
\begin{center}
\begin{tabular}{|c|cccccc|}
  \hline
   grid & $n^{iter}_{\text{\tiny{MINRES}}}$ & $t^{\text{sec}}$ & $\|{\mathcal{R}_1}_1\|_{\Omega}$
   & $\|{\mathcal{R}_2}_1\|_{\Omega}$ & $\mathcal{M}^{\oplus_1}_{|\cdot|}$ & $I_{\text{eff}}^1$  \\
  \hline
   $9 \times 9$     & 14 & 0.00  & 1.238e+02 & 3.444e+00 & 4.426e+01 & 1.908 \\ 
   $27 \times 27$   & 12 & 0.02  & 4.566e+01 & 7.376e-02 & 1.464e+01 & 1.550 \\ 
   $81 \times 81$   & 10 & 0.21  & 1.540e+01 & 5.560e-03 & 4.910e+00 & 1.485 \\ 
   $243 \times 243$ & 10 & 2.08  & 5.140e+00 & 3.659e-04 & 1.637e+00 & 1.466 \\ 
   $729 \times 729$ & 8  & 15.84 & 1.714e+00 & 2.295e-05 & 5.455e-01 & 1.460 \\ 
  \hline
\end{tabular}
\end{center}
\caption{Majorant $\mathcal{M}^{\oplus_1}_{|\cdot|}$ and its parts
(Example 2).}
\label{tab:Ex2:1}
\end{table}
\begin{table}[!ht]
\begin{center}
\begin{tabular}{|c|cccccc|}
  \hline
   grid & $n^{iter}_{\text{\tiny{MINRES}}}$ & $t^{\text{sec}}$ & $\|{\mathcal{R}_1}_2\|_{\Omega}$
   & $\|{\mathcal{R}_2}_2\|_{\Omega}$ & $\mathcal{M}^{\oplus_2}_{|\cdot|}$ & $I_{\text{eff}}^2$  \\
  \hline
   $9 \times 9$     & 9 & 0.00  & 7.953e+01 & 2.209e+00 & 2.844e+01 & 1.880 \\ 
   $27 \times 27$   & 9 & 0.02  & 2.930e+01 & 4.737e-02 & 9.394e+00 & 1.523 \\ 
   $81 \times 81$   & 9 & 0.19  & 9.883e+00 & 3.555e-03 & 3.151e+00 & 1.460 \\ 
   $243 \times 243$ & 8 & 1.73  & 3.299e+00 & 2.339e-04 & 1.050e+00 & 1.441 \\ 
   $729 \times 729$ & 8 & 15.64 & 1.100e+00 & 1.467e-05 & 3.501e-01 & 1.435 \\ 
  \hline
\end{tabular}
\end{center}
\caption{Majorant $\mathcal{M}^{\oplus_2}_{|\cdot|}$ and its parts
(Example 2).}
\label{tab:Ex2:2}
\end{table}
\begin{table}[!ht]
\begin{center}
\begin{tabular}{|c|cccccc|}
  \hline
   grid & $n^{iter}_{\text{\tiny{MINRES}}}$ & $t^{\text{sec}}$ & $\|{\mathcal{R}_1}_3\|_{\Omega}$
   & $\|{\mathcal{R}_2}_3\|_{\Omega}$ & $\mathcal{M}^{\oplus_3}_{|\cdot|}$ & $I_{\text{eff}}^3$  \\
  \hline
   $9 \times 9$     & 8 & 0.00  & 4.613e+01 & 1.277e+00 & 1.649e+01 & 1.905 \\ 
   $27 \times 27$   & 8 & 0.02  & 1.696e+01 & 2.745e-02 & 5.437e+00 & 1.541 \\ 
   $81 \times 81$   & 7 & 0.15  & 5.719e+00 & 2.046e-03 & 1.823e+00 & 1.477 \\ 
   $243 \times 243$ & 7 & 1.56  & 1.909e+00 & 1.345e-04 & 6.079e-01 & 1.457 \\ 
   $729 \times 729$ & 6 & 12.34 & 6.364e-01 & 8.436e-06 & 2.026e-01 & 1.451 \\ 
  \hline
\end{tabular}
\end{center}
\caption{Majorant $\mathcal{M}^{\oplus_3}_{|\cdot|}$ and its parts
(Example 2).}
\label{tab:Ex2:3}
\end{table}
\begin{table}[!ht]
\begin{center}
\begin{tabular}{|c|cccccc|}
  \hline
   grid & $n^{iter}_{\text{\tiny{MINRES}}}$ & $t^{\text{sec}}$ & $\|{\mathcal{R}_1}_4\|_{\Omega}$
   & $\|{\mathcal{R}_2}_4\|_{\Omega}$ & $\mathcal{M}^{\oplus_4}_{|\cdot|}$ & $I_{\text{eff}}^4$  \\
  \hline
   $9 \times 9$     & 10 & 0.00  & 1.624e+01 & 4.474e-01 & 5.801e+00 & 1.958 \\ 
   $27 \times 27$   & 9  & 0.02  & 5.950e+00 & 9.645e-03 & 1.908e+00 & 1.582 \\ 
   $81 \times 81$   & 9  & 0.19  & 2.007e+00 & 7.120e-04 & 6.398e-01 & 1.516 \\ 
   $243 \times 243$ & 9  & 1.90  & 6.698e-01 & 4.680e-05 & 2.133e-01 & 1.496 \\ 
   $729 \times 729$ & 8  & 15.88 & 2.233e-01 & 2.934e-06 & 7.108e-02 & 1.490 \\ 
  \hline
\end{tabular}
\end{center}
\caption{Majorant $\mathcal{M}^{\oplus_4}_{|\cdot|}$ and its parts
(Example 2).}
\label{tab:Ex2:4}
\end{table}
\begin{table}[!ht]
\begin{center}
\begin{tabular}{|c|cccccc|}
  \hline
   grid & $n^{iter}_{\text{\tiny{MINRES}}}$ & $t^{\text{sec}}$ & $\|{\mathcal{R}_1}_5\|_{\Omega}$
   & $\|{\mathcal{R}_2}_5\|_{\Omega}$ & $\mathcal{M}^{\oplus_5}_{|\cdot|}$ & $I_{\text{eff}}^5$  \\
  \hline
   $9 \times 9$     & 9 & 0.00  & 5.878e+00 & 1.611e-01 & 2.099e+00 & 1.970 \\ 
   $27 \times 27$   & 9 & 0.02  & 2.146e+00 & 3.485e-03 & 6.880e-01 & 1.586 \\ 
   $81 \times 81$   & 7 & 0.15  & 7.236e-01 & 2.542e-04 & 2.307e-01 & 1.520 \\ 
   $243 \times 243$ & 7 & 1.54  & 2.415e-01 & 1.669e-05 & 7.690e-02 & 1.500 \\ 
   $729 \times 729$ & 7 & 14.17 & 8.052e-02 & 1.046e-06 & 2.563e-02 & 1.493 \\ 
  \hline
\end{tabular}
\end{center}
\caption{Majorant $\mathcal{M}^{\oplus_5}_{|\cdot|}$ and its parts
(Example 2).}
\label{tab:Ex2:5}
\end{table}
\begin{table}[!ht]
\begin{center}
\begin{tabular}{|c|cccccc|}
  \hline
   grid & $n^{iter}_{\text{\tiny{MINRES}}}$ & $t^{\text{sec}}$ & $\|{\mathcal{R}_1}_6\|_{\Omega}$
   & $\|{\mathcal{R}_2}_6\|_{\Omega}$ & $\mathcal{M}^{\oplus_6}_{|\cdot|}$ & $I_{\text{eff}}^6$  \\
  \hline
   $9 \times 9$     & 11 & 0.00  & 2.621e+00 & 7.132e-02 & 9.351e-01 & 1.991 \\ 
   $27 \times 27$   & 10 & 0.02  & 9.522e-01 & 1.550e-03 & 3.053e-01 & 1.597 \\ 
   $81 \times 81$   & 9  & 0.19  & 3.211e-01 & 1.114e-04 & 1.024e-01 & 1.529 \\ 
   $243 \times 243$ & 8  & 1.73  & 1.072e-01 & 7.312e-06 & 3.412e-02 & 1.509 \\ 
   $729 \times 729$ & 8  & 15.81 & 3.573e-02 & 4.583e-07 & 1.137e-02 & 1.503 \\ 
  \hline
\end{tabular}
\end{center}
\caption{Majorant $\mathcal{M}^{\oplus_6}_{|\cdot|}$ and its parts
(Example 2).}
\label{tab:Ex2:6}
\end{table}
\begin{table}[!ht]
\begin{center}
\begin{tabular}{|c|cccccc|}
  \hline
   grid & $n^{iter}_{\text{\tiny{MINRES}}}$ & $t^{\text{sec}}$ & $\|{\mathcal{R}_1}_7\|_{\Omega}$
   & $\|{\mathcal{R}_2}_7\|_{\Omega}$ & $\mathcal{M}^{\oplus_7}_{|\cdot|}$ & $I_{\text{eff}}^7$  \\
  \hline
   $9 \times 9$     & 12 & 0.00  & 1.359e+00 & 3.670e-02 & 4.846e-01 & 2.023 \\ 
   $27 \times 27$   & 11 & 0.02  & 4.913e-01 & 8.015e-04 & 1.575e-01 & 1.615 \\ 
   $81 \times 81$   & 9  & 0.19  & 1.656e-01 & 5.669e-05 & 5.280e-02 & 1.546 \\ 
   $243 \times 243$ & 9  & 1.89  & 5.528e-02 & 3.716e-06 & 1.760e-02 & 1.526 \\ 
   $729 \times 729$ & 8  & 15.88 & 1.843e-02 & 2.329e-07 & 5.867e-03 & 1.520 \\ 
  \hline
\end{tabular}
\end{center}
\caption{Majorant $\mathcal{M}^{\oplus_7}_{|\cdot|}$ and its parts
(Example 2).}
\label{tab:Ex2:7}
\end{table}
\begin{table}[!ht]
\begin{center}
\begin{tabular}{|c|cccccc|}
  \hline
   grid & $n^{iter}_{\text{\tiny{MINRES}}}$ & $t^{\text{sec}}$ & $\|{\mathcal{R}_1}_8\|_{\Omega}$
   & $\|{\mathcal{R}_2}_8\|_{\Omega}$ & $\mathcal{M}^{\oplus_8}_{|\cdot|}$ & $I_{\text{eff}}^8$  \\
  \hline
   $9 \times 9$     & 12 & 0.00  & 7.839e-01 & 2.098e-02 & 2.792e-01 & 2.064 \\ 
   $27 \times 27$   & 11 & 0.02  & 2.816e-01 & 4.606e-04 & 9.028e-02 & 1.638 \\ 
   $81 \times 81$   & 11 & 0.23  & 9.492e-02 & 3.200e-05 & 3.026e-02 & 1.569 \\ 
   $243 \times 243$ & 10 & 2.05  & 3.168e-02 & 2.095e-06 & 1.009e-02 & 1.548 \\ 
   $729 \times 729$ & 10 & 19.27 & 1.056e-02 & 1.312e-07 & 3.362e-03 & 1.542 \\ 
  \hline
\end{tabular}
\end{center}
\caption{Majorant $\mathcal{M}^{\oplus_8}_{|\cdot|}$ and its parts
(Example 2).}
\label{tab:Ex2:8}
\end{table}
\begin{table}[!ht]
\begin{center}
\begin{tabular}{|c|cccccc|}
  \hline
   & $n^{iter}_{\text{\tiny{MINRES}}}$ & $t^{\text{sec}}$ & $\|{\mathcal{R}_1}\|$
   & $\|{\mathcal{R}_2}\|$ & $\mathcal{M}^{\oplus}_{|\cdot|}$ & $I_{\text{eff}}$  \\
  \hline
   $k=0$   & -  & -     & 8.821e-01 & 1.183e-05 & 1.986e-01 & 1.002 \\ 
   $k=1$   & 8  & 15.84 & 1.714e+00 & 2.295e-05 & 5.455e-01 & 1.460 \\ 
   $k=2$   & 8  & 15.64 & 1.100e+00 & 1.467e-05 & 3.501e-01 & 1.435 \\ 
   $k=3$   & 6  & 12.34 & 6.364e-01 & 8.436e-06 & 2.026e-01 & 1.451 \\ 
   $k=4$   & 8  & 15.88 & 2.233e-01 & 2.934e-06 & 7.108e-02 & 1.490 \\ 
   $k=5$   & 7  & 14.17 & 8.052e-02 & 1.046e-06 & 2.563e-02 & 1.493 \\ 
   $k=6$   & 8  & 15.81 & 3.573e-02 & 4.583e-07 & 1.137e-02 & 1.503 \\ 
   $k=7$   & 8  & 15.88 & 1.843e-02 & 2.329e-07 & 5.867e-03 & 1.520 \\ 
   $k=8$   & 10 & 19.27 & 1.056e-02 & 1.312e-07 & 3.362e-03 & 1.542 \\ 
  \hline
   overall & -  & -     & 4.403e+00 & 5.886e-05 & 1.402e+00 & 1.404 \\ 
  \hline
\end{tabular}
\end{center}
\caption{The overall majorant $\mathcal{M}^{\oplus}_{|\cdot|}$ and its parts
computed on a $729 \times 729$-mesh
(Example 2).}
\label{tab:Ex2:global}
\end{table}

\section*{Acknowledgment}

The research was supported by the Austrian Science Fund (FWF) under the grant W1214-N15,
project DK4, as well as by the strategic program ``Innovatives O\"O 2010 plus'' by the Upper
Austrian Government, and by the Austrian Academy of Sciences.

\bibliographystyle{abbrv}
\bibliography{LRW1arxiv}
\end{document}